\documentclass[hyperref]{article}
\usepackage{graphicx} %% Package for Figure
\usepackage{float} %% Package for Float
\usepackage{amssymb}
\usepackage{amsmath}
\usepackage{mathtools}
\usepackage{bm} %% Bold Mathematical Symbols
\usepackage[colorlinks,linkcolor=cyan,citecolor=cyan]{hyperref}
\usepackage{extarrows}
\usepackage[hang,flushmargin]{footmisc} %% Let the footnote not indentation
\usepackage[square,comma,sort&compress,numbers]{natbib} %% Sort of References
\usepackage{mathrsfs} %% Swash letter
\usepackage[font=footnotesize,skip=0pt,textfont=rm,labelfont=rm]{caption,subcaption}
\usepackage{booktabs} %% tables with three lines
\usepackage{tocloft}
\usepackage[all]{xy}
\usepackage{indentfirst}
\usepackage{tikz-cd}
\usepackage{graphicx}
\usepackage{color}
\usepackage{ulem}

\makeatletter

\newcommand{\Rmnum}[1]{\expandafter\@slowromancap\romannumeral #1@}
\makeatother

\usepackage{amsthm}

\newtheorem{Thm}{Theorem}[section]
\newtheorem{Cor}[Thm]{Corollary}
\newtheorem{Prop}[Thm]{Proposition}
\newtheorem{Lem}[Thm]{Lemma}
\newtheorem{Assumption}[Thm]{Assumption}
\newtheorem{Notation}[Thm]{Notation}

\theoremstyle{definition}
\newtheorem{Def}[Thm]{Definition}
\newtheorem{Rem}[Thm]{Remark}
\newtheorem{Ex1}[Thm]{Example}

\DeclareMathOperator{\HH}{HH}

\DeclareMathOperator{\rank}{rank}

\DeclareMathOperator{\Sp}{Sp}
\DeclareMathOperator{\Ker}{Ker}
\renewcommand{\Im}{\mathrm{Im}}
\renewcommand{\sl}{\mathfrak{sl}}
\renewcommand{\dim}{\mathrm{dim}}

\newcommand{\pirank}{\pi_1\!\text{-}\!\rank}

\newcommand{\gl}{\mathfrak{gl}}

\graphicspath{{figure/}}                                 %% Path of Figures
\usepackage[a4paper]{geometry}                           %% Paper size
\begin{document}
\title{\bf The Hochschild cohomology groups under gluing arrows}
\author{Yuming Liu$^a$, Lleonard Rubio y Degrassi$^b$, Can Wen$^{a,*}$}
\maketitle

$^a$ Yuming Liu and Can Wen, School of Mathematical Sciences, Laboratory of Mathematics and Complex Systems,
Beijing Normal University, Beijing 100875, P. R. China.

$^b$ Lleonard Rubio y Degrassi, Department of Mathematics, Uppsala University, Box 480, 75106, Uppsala, Sweden.

E-mail addresses: ymliu@bnu.edu.cn (Y. Liu); lleonard.rubio@math.uu.se (L. Rubio y Degrassi); \\ cwen@mail.bnu.edu.cn (C. Wen).

$^*$ Corresponding author.

%MS+++++++++++++++++++++ Abstract +++++++++++++++++++++++++

\renewcommand{\thefootnote}{\alph{footnote}}
\setcounter{footnote}{-1} \footnote{\it{Mathematics Subject
		Classification(2020)}: 16E40, 16G10.}
\renewcommand{\thefootnote}{\alph{footnote}}
\setcounter{footnote}{-1} \footnote{ \it{Keywords}: Fundamental group; Gluing arrows; Hochschild cohomology; Monomial algebra.}

{\noindent\small{\bf Abstract:} In a previous paper \cite{LRW} we have compared the Hochschild cohomology groups of finite dimensional monomial algebras under gluing two idempotents. In the present paper, we compare the Hochschild cohomology groups of finite dimensional monomial algebras under gluing two arrows.}

%\vspace{1ex}
%{\noindent\small{\bf Keywords:}
%    Keywords1; Keywords2;...}
%MS++++++++++++++++++++++++++++++ Main body ++++++++++++++++++++

\section{Introduction}
In a previous paper \cite{LRW}, we have studied the behaviour of Hochschild cohomology groups and of the fundamental groups of finite dimensional monomial algebras under gluing two idempotents. Recall that if $A=kQ_A/I_A$ is a finite dimensional algebra and $B\simeq kQ_B/I_B$ is a subalgebra of $A$ obtained by gluing two arbitrary idempotents of $A$ (that is, gluing two arbitrary vertices of $Q_A$), then the canonical embedding $\phi: B \to A$ preserves the Jacobson radicals of algebras, in other words, $\phi(\mathrm{rad}\, B)=\mathrm{rad}\, A$. Suppose that $A$ (hence also $B$) is a monomial algebra. In \cite{LRW}, we provided explicit formulas between dimensions of $\HH^i(A)$ and $\HH^i(B)$ for $i=0,1$ in terms of some combinatorial datum; when gluing a source vertex and a sink vertex, we also described the relationship between Lie algebra structures of $\HH^1(A)$ and $\HH^1(B)$.
	
	The main aim of the present paper is to study the behaviours of Hochschild cohomology groups and of the fundamental
	groups when we glue two arrows from a finite dimensional monomial algebra $A=kQ_A/I_A$.
	
The idea of gluing arrows is similar to that of gluing idempotents and has appeared in Guo's master thesis \cite{G} when he generalizes some results on stable equivalences modulo projective modules (or modulo semisimple modules) in \cite{KL} to stable equivalences modulo more general modules. If $A=kQ_A/I_A$ is a finite dimensional algebra, then by gluing two arrows of $Q_A$ we can form a subalgebra (having the same identity element as $A$) $B$ of $A$. However, in this case the canonical embedding $\phi: B \to A$ does not preserve the Jacobson radicals anymore; it turns out that $\phi$ preserves the square of the Jacobson radicals. Another important difference between these two kinds of gluings is that gluing idempotents produce new relations given by paths of length two while gluing arrows will produce new relations given by paths of length two or three. Nonetheless, it is rather interesting that the most results in \cite{LRW} for gluing vertices can be generalized to the situation of gluing arrows.
	
This paper is organized as follows: in Section \ref{prelim} we introduce some notation that will be used for the rest of the paper. We also give some background.
	 In Section \ref{firstHoch} we study the behaviour of the first Hochschild cohomology in case of gluing a source and a sink arrow. This culminates in Theorem \ref{hh1-gluing-source-sink-arrows}
	  in which we compare the Lie algebra structures. 		
	 Although in case of gluing of two arbitrary arrows structural results do not hold, cf. Remark \ref{dim-Image}, in  Theorem \ref{hh1-glue-two-arrows}  we prove that  the same dimension formula
	 holds. In addition, if the characteristic of the field is zero, then we give some applications of the main results in Corollaries \ref{one-dim-summand-source-sink} and \ref{one-dim-summand-rad-aquare-zero}. In
	 Section \ref{centergluingarrows} we study how gluing changes the center. In Section \ref{fundsec} we study the relation between gluing arrows and the $\pi_1$-ranks.
	  In particular, we will establish a connection between the dual fundamental group and the first Hochschild cohomology in case of gluing a source and a sink arrow. In Section \ref{higherhoch} we consider how higher Hochschild cohomology groups change with respect to the gluing of a source and a sink arrow.	
	  In Section \ref{Ex} we give various examples to illustrate our definitions and results.

\section{Preliminaries}\label{prelim}
Throughout this paper, we follow the same notation of our previous work \cite{LRW}. Let us recall some details: let $A$ be a finite dimensional algebra of the form $kQ_A/I_A$, where $k$ is a field with char$(k)\geq 0$, $Q_A$ is a finite quiver (with vertex set $(Q_A)_0$ and arrow set $(Q_A)_1$) and $I_A$ is an admissible ideal in the path algebra $kQ_A$. Denote the vertices of $Q_A$ by $e_1,\dots,e_n$ without distinguishing the notation of the idempotents of $A$, and denote the arrows of $Q_A$ by the Greek letters $\alpha,\beta,\gamma,\dots$. For an arrow $\alpha$, $s(\alpha)$ (resp. $t(\alpha)$) denotes the starting vertex (resp. the ending vertex) of $\alpha$. For a path $p=\alpha_n\dots \alpha_1$ in $Q_A$, $s(p):=s(\alpha_1)$ (resp. $t(p):=t(\alpha_n)$) denotes the starting vertex (resp. the ending vertex) of $p$. Two paths $\epsilon,\gamma$ of $Q_A$ are called parallel if $s(\epsilon)=s(\gamma)$ and $t(\epsilon)=t(\gamma)$. We denote the pair $(\epsilon, \gamma)$ by  $\epsilon \| \gamma$. If $\epsilon$ and $\gamma$ are not parallel, then we denote by $\epsilon \nparallel \gamma$. If $X,Y$ are sets of paths of $Q_A$, then we denote by $X \| Y$ the set of parallel paths consisting of the couples $\epsilon \| \gamma$ with $\epsilon \in X$ and $\gamma \in Y$, and denote by $k(X \| Y)$ the $k$-vector space with basis $X \| Y$. We denote by $\mathrm{dim}\,V$ the dimension of a $k$-vector space $V$.

{\bf Gluing arrows.} We now fix the following notations. Let us assume that the quiver $Q_A$ has $n$ $(n\ge 4)$ vertices (denoted by $e_1,\dots,e_n$) and $m+2$ $(m\geq 0)$ arrows (denoted by $\alpha,\beta,\gamma_1,\dots,\gamma_m$), where $\alpha: e_1 \to e_2$ and $\beta: e_{n-1} \to e_n$. We will throughout assume that the four vertices of $\alpha$ and $\beta$ are pairwise different (cf. Remark \ref{exculde-cases}). Let $B$ be an algebra obtained from $A$ by gluing $\alpha$ and $\beta$, that is, $B$ is identified as a subalgebra of $A$
generated by $f_1:=e_1+e_{n-1},f_2=e_2+e_n,f_i=e_i$ for $3 \le i \le n-2$ and $\gamma^*=\alpha +\beta,\gamma_j^*=\gamma_j$ for $1 \le j \le m$. Then it is not hard to verify that the subalgebra $B$ has the same identity element as $A$ (that is, $1_A=e_1+e_2\cdots+e_n=f_1+f_2+\cdots+f_{n-2}=1_B$). In addition, $B$ is isomorphic to $kQ_B/I_B$, where $Q_B$ is the quiver obtained from $Q_A$ by identifying the arrows $\alpha$ and $\beta$ (thus also identifying the vertices $e_1$ and $e_{n-1}$, and identifying the vertices $e_2$ and $e_n$) and where $I_B$ is an admissible ideal generated by the elements in $I_A$ and all newly formed paths of length 2 through $f_1$ or $f_2$ and all newly formed paths of length 3 through $\gamma^*$. More specifically, the ideal $I_B$ is generated by $I_A$ and  $Z_{new}$, where $Z_{new}$ can be described as follows for $\eta,\mu,\lambda,\xi \in (Q_A)_1$:
\begin{equation*}
\begin{split}
Z_{new}&=\{\lambda^*\eta^*,\xi^*\mu^*,\xi^*\gamma^*\eta^*\ |\ t(\eta)=e_1,t(\mu)=e_2(\mu\neq \alpha),s(\lambda)=e_{n-1}(\lambda\neq \beta),s(\xi)=e_n\}\\ &\cup\{\lambda^*\eta^*,\xi^*\mu^*,\xi^*\gamma^*\eta^*\ |\ t(\eta)=e_{n-1},t(\mu)=e_n(\mu\neq \beta),s(\lambda)=e_1(\lambda \neq \alpha),s(\xi)=e_2\}.
\end{split}
\end{equation*}

\iffalse
\begin{center}
	\begin{tikzcd}
	& {} & {} & {} & {} &&& {} & {} & {} \\
	{Q_A:}\hspace{-3em} & {e_1\bullet} & {e_2\bullet} & {e_{n-1}\bullet} & {\bullet e_n} && {Q_B:}\hspace{-3em} & {} & {f_1\bullet} & {\bullet f_2} & {} \\
	&&& {} &&&&& {}
	\arrow["\alpha", from=2-2, to=2-3]
	\arrow["\beta", from=2-4, to=2-5]
	\arrow["\cdots"{description}, draw=none, from=2-3, to=2-4]
	\arrow["\xi"', shift left=1, from=2-5, to=1-5]
	\arrow["\mu", shift left=2, from=1-3, to=2-3]
	\arrow["{\eta^*}"', from=1-9, to=2-9]
	\arrow["{\xi^*}", from=2-10, to=2-11]
	\arrow["{\gamma^*}", shift left=2, from=2-9, to=2-10]
	\arrow["\cdots"{description}, shift right=2, draw=none, from=2-9, to=2-10]
	\arrow["{\lambda^*}"', from=2-9, to=2-8]
	\arrow["{\mu^*}", from=1-10, to=2-10]
	\arrow["\eta"', from=1-2, to=2-2]
	\arrow["\lambda", from=2-4, to=1-4]
	\end{tikzcd}$(*)$
\end{center}
\fi
Note that when gluing two arrows under the above assumption, the length of each path in $Z_{new}$ is 2 or 3; moreover, we have $\dim\,B=\dim\,A-3$ and rad$^2 B$=rad$^2 A$ (cf. \cite[Section 2.4]{G}). \noindent  It is worth to mention that, from an inverse gluing operation point of view (cf. \cite[Example 3.32]{LRW}), by inverse gluing idempotents operation, we can reduce a relation in $I_B$ given by a path of length $2$, by inverse gluing arrows operation, we can also reduce a relation in $I_B$ given by a path of length $3$ through a node arrow (cf. Definition \ref{node arrow}).

\begin{Ex1} The algebra $B$ is obtained from a radical cube zero algebra $A$ by gluing $\alpha$ and $\beta$:
	\begin{center}
		\begin{tikzcd}
		Q_A:\hspace{-2em} & e_1\bullet \arrow[r, "\mu" description, shift left=4] \arrow[r, "\alpha" description] & e_2\bullet \arrow[l, "\eta" description, shift left=4] & e_3\bullet \arrow[l, "\lambda"'] \arrow[r, "\beta" description] \arrow[r, "b" description, shift left=4] & \bullet e_4 \arrow[l, "\xi" description, shift left=4] &  & Q_B:\hspace{-2em} & f_1\bullet \arrow[r, "\gamma^*" description, shift left=2] \arrow[r, "\mu^*" description, shift left=6] \arrow[r, "b^*" description, shift left=10] \arrow[r, "\lambda^*" description, shift right=10] & \bullet f_2 \arrow[l, "\eta^*" description, shift left=2] \arrow[l, "\xi^*" description, shift left=6]
		\end{tikzcd},
	\end{center}
where $Z_{new}=\{\lambda^*\eta^*,b^*\eta^*,\xi^*\lambda^*,\xi^*\mu^*,\xi^*\gamma^*\eta^* \} \cup \{\mu^*\xi^*,\eta^*b^*,\eta^*\gamma^*\xi^* \}$.
\end{Ex1}

\begin{Ex1}\label{ss-same-block}
The algebra $B$ is obtained from a path algebra $A=kQ_A$ by gluing $\alpha$ and $\beta$:
	\begin{center}
		\begin{tikzcd}
		Q_A:\hspace{-2em} & e_1\bullet \arrow[r, "\alpha"] & e_2\bullet \arrow[r, "\eta"] & e_3\bullet \arrow[r, "\beta"] & \bullet e_4 &  & Q_B:\hspace{-2em} & f_1\bullet \arrow[r, "\gamma^*", shift left=2] & \bullet f_2 \arrow[l, "\eta^*", shift left=2]
		\end{tikzcd},
	\end{center}
where $Z_{new}=\{\eta^*\gamma^*\eta^* \}$.
\end{Ex1}

\begin{Ex1}\label{ss-different-blocks}
The algebra $B$ is obtained from a path algebra $A=kQ_A$ by gluing $\alpha$ and $\beta$:
	\begin{center}
		\begin{tikzcd}
		Q_A:\hspace{-2em} & e_4\bullet \arrow[r, "\delta"] & e_5\bullet \arrow[r, "\beta"]& e_6\bullet & e_1\bullet \arrow[r, "\alpha"] & e_2\bullet \arrow[r, "\varepsilon"]& e_3 \bullet \\ Q_B:\hspace{-2em} & f_4\bullet \arrow[r, "\delta^*"] & f_1\bullet \arrow[r, "\gamma^*"] & f_2\bullet \arrow[r, "\varepsilon^*"] & \bullet f_3
		\end{tikzcd},
	\end{center}
where $Z_{new}=\{\varepsilon^*\gamma^*\delta^* \}$. Note that in this case the quiver $Q_A$ has two connected components and therefore the algebra $A$ has two blocks.
\end{Ex1}

\begin{Rem} \label{quiver-morphism} Similar to the gluing idempotents situation in \cite[Section 2]{LRW}, we have the following elementary observation: The gluing arrows operation induces a quiver morphism $\varphi:Q_{A} \to Q_{B}$ which is given by the following formula: $\varphi(e_{i})=f_{i}$ for $3\le i \le n-2$, $\varphi(e_{1})=\varphi(e_{n-1})=f_{1}$, $\varphi(e_2)=\varphi(e_n)=f_2$, and $\varphi(\alpha)=\varphi(\beta)=\gamma^*$ and $\varphi(\gamma_i)=\gamma_i^*$ for each $\gamma_i \in (Q_A)_1\backslash \{\alpha,\beta\}$.
\end{Rem}

\begin{Rem} \label{exculde-cases}
Since we expect that the subalgebra $B$ is a unital algebra and has the same identity element as $A$ and the quiver $Q_B$ is obtained from $Q_A$ by gluing the arrows $\alpha$ and $\beta$, it is necessary to assume that the four vertices of $\alpha$ and $\beta$ are pairwise different. In particular, we will exclude the following situations.
\begin{itemize}	
\item[$(1)$] $\alpha$ or $\beta$ is a loop. For example, let $A$ be defined by the following quiver
\begin{center}
		\begin{tikzcd}
		Q_A:\hspace{-2em} & e_1\bullet \arrow[r, "\alpha"] & \bullet e_2 \arrow["\beta", loop, distance=2em, in=55, out=125]
		\end{tikzcd}
	\end{center}
with relation $\beta^2=0$. It is easy to verify that the subalgebra $B$ (generated by $e_1+e_2$, $e_2+e_2$ and $\alpha+\beta$) and $A$ (hence $Q_B$ and $Q_A$) are equal when $\mathrm{char}(k)=0$. Note that if we define the subalgebra $B$ generated by $e_1+e_2$ and $\alpha+\beta$, then $B$ has the same identity element as $A$ and $Q_B$ is obtained from $Q_A$ by gluing the arrows $\alpha$ and $\beta$ and $B$ is radical cube zero; however, the ideal $I_B$ is not generated by $I_A$ and $Z_{new}$.
\item[$(2)$] Both $\alpha$ and $\beta$ are not loops but they share one common vertex. For example, let $A$ be a path algebra defined by the following quiver
	\begin{center}
		\begin{tikzcd}
		 & & e_2\bullet                                        &                        \\
		Q_A:\hspace{-2em} &e_4\bullet \arrow[r, "\eta"] & e_1\bullet \arrow[u, "\alpha"] \arrow[r, "\beta"] & \bullet e_3
		\end{tikzcd}
	\end{center}
with relation $\alpha\eta=0$. If $\mathrm{char}(k)=2$, then $e_1+e_1=0$. It follows easily that the subalgebra $B$ of $A$ generated by $e_1+e_1$, $e_2+e_3, e_4, \eta$ and $\alpha+\beta$ is a non-unital algebra. Note that if we define the subalgebra $B$ generated by $e_1$, $e_2+e_3$, $e_4$, $\eta$ and $\alpha+\beta$, then $B$ has the same identity element as $A$ and $Q_B$ is obtained from $Q_A$ by gluing the arrows $\alpha$ and $\beta$; however, the ideal $I_B$ is equal to zero and clearly not generated by $I_A$ and $Z_{new}$.
\end{itemize}	
\end{Rem}

We are particularly interested in the case of gluing a source arrow and a sink arrow.

\begin{Def} {\rm(\cite[Section 4.1]{G})} \label{source-sink arrow}
	An arrow $\alpha$ in $Q_A$ is called as a source arrow if the following three conditions are satisfied:
	
	(1) s($\alpha$) is a source vertex;
	
	(2) except $\alpha$, there is neither an arrow starting from s($\alpha$) nor an arrow ending at t($\alpha$).
	
	Dually, an arrow $\beta$ in $Q_A$ is called as a sink arrow if the following three conditions are satisfied:
	
	(1) t($\beta$) is a sink vertex;
	
	(2) except $\beta$, there is neither an arrow starting from s($\beta$) nor an arrow ending at t($\beta$).
	
\end{Def}

\begin{Def} \rm{(\cite[Section 4.1]{G})} \label{node arrow}
	An arrow $\gamma$ in $Q_A$ is called as a node arrow if the following three conditions are satisfied:
	
	(1) $\gamma$ is not a source arrow;
	
	(2) $\gamma$ is not a sink arrow;
	
	(3) all paths of length $3$ in $Q_A$ of the form $\cdot\to\cdot\stackrel{\gamma}{\longrightarrow}\cdot\to\cdot$ belong to $I_A$.
\end{Def}

\begin{Rem}
\begin{itemize}	
\item[$(1)$] Let $\alpha$ be a source arrow (but not a sink arrow) and $\beta$ be a sink arrow (but not a source arrow) in $Q_A$. If $\gamma^*$ is the arrow obtained by gluing $\alpha$ and $\beta$, then $\gamma^*$ is a node arrow in $Q_B$.
\item[$(2)$] If we glue a source arrow and a sink arrow, then either $Z_{new}$ is empty or the length of each path in $Z_{new}$ is just $3$ (cf. Examples \ref{ss-same-block} and \ref{ss-different-blocks}).
\end{itemize}
\end{Rem}

The following lemma is similar to \cite[Lemma 3.1-3.2]{LRW} and its proof is straightforward.

\begin{Lem} \label{parallel-arrows}
Let $B$ be a subalgebra of a finite dimensional algebra $A$ obtained by gluing two arrows $\alpha$ and $\beta$ of $Q_A$. For $\eta,\delta\in (Q_A)_1$, the following statements hold. \begin{itemize}	
\item[$(1)$] If $\eta\|\delta$ in $Q_A$, then $\eta^{*} \| \delta^{*}$ in $Q_B$.
\item[$(2)$] Suppose that $\alpha$ is a source arrow and $\beta$ is a sink arrow. Then $\eta \| \delta$ if and only if $\eta^{*} \| \delta^{*}$, except that $\eta=\alpha, \delta=\beta$ or $\eta=\beta, \delta=\alpha$.
\end{itemize}
\end{Lem}

{\bf Hochschild cochain complex of monomial algebras.}
Let $A$ be a finite dimensional monomial $k$-algebra, that is, a finite dimensional $k$-algebra which is isomorphic to a quotient $kQ_A/I$ of a path algebra where the two-sided ideal $I$ of $kQ_A$ is generated by a set $Z=Z_A$ of paths of length $\ge$ 2. We shall assume that $Z$ is minimal, that is, no proper subpath of a path in $Z$ is again in $Z$. Let $\mathcal{B}=\mathcal{B}_{A}$ be the set of paths of $Q_A$ which do not contain any element of $Z$ as a subpath. It is clear that the (classes modulo $I$ of) elements of $\mathcal{B}$ form a basis of $A$.

Recall that Strametz \cite{STR} gave a method to compute the Hochschild cohomology groups $\mathrm{HH}^{n}(A)$ in degree $n=0,1$ of a monomial algebra $A$  using parallel paths in $Q:=Q_A$.

\begin{Prop}\label{Strametz-proposition2.6} {\rm(\cite[Proposition 2.6]{STR})}
\label{straCmon}
	Let $A$ be a finite dimensional monomial algebra. Consider the following cochain complex (denoted by $\mathcal{C}_{mon}$):
	\begin{center}
		\begin{tikzcd}
		0 \arrow[r] & k(Q_{0}\| \mathcal{B}) \arrow[r, "\delta^{0}"] & k(Q_{1}\| \mathcal{B}) \arrow[r, "\delta^{1}"] & k(Z\| \mathcal{B}) \arrow[r, "\delta^{2}"] & \cdots,
		\end{tikzcd}	
	\end{center}
	where the differentials are given by
	\begin{equation*}
	\begin{split}
	\delta^{0}:k(Q_{0}\| \mathcal{B}) &\to k(Q_{1}\| \mathcal{B})\\
	e\| \gamma &\mapsto \sum_{s(a)=e,a\gamma \in \mathcal{B}}a\| a\gamma - \sum_{t(a)=e,\gamma a \in \mathcal{B}}a\| \gamma a,\\
	\delta^{1}:k(Q_{1}\| \mathcal{B}) &\to k(Z \| \mathcal{B})\\
	a\| \gamma &\mapsto \sum_{r \in Z}r\| r^{a \| \gamma},
	\end{split}
	\end{equation*}
where $r^{a \| \gamma}$ denotes the sum of all paths in $\mathcal{B}$ obtained by replacing each appearance of the arrow $a$ in $r$ by the path $\gamma$. Then we have $\mathrm{HH}^{0}(A)\simeq\mathrm{Ker}(\delta^{0})$ and $\mathrm{HH}^{1}(A)\simeq\mathrm{Ker}(\delta^{1})/ \mathrm{Im}(\delta^{0})$.
\end{Prop}

\begin{Thm}\label{Strametz-theorem2.7} {\rm(\cite{STR})}
Let $A$ be a finite dimensional monomial algebra. Then the bracket $[\ ,\ ]: k(Q_1\|\mathcal{B})\times k(Q_1\|\mathcal{B})\to k(Q_1\|\mathcal{B})$ given by $$\left[a \| \gamma,b \| \epsilon \right] = b\| \epsilon^{a \| \gamma}-a\| \gamma^{b \| \epsilon}\quad \quad (a\| \gamma,  b\| \epsilon \in Q_{1} \| \mathcal{B})$$
induces a Lie algebra structure on $\mathrm{HH}^{1}(A)\simeq\mathrm{Ker}(\delta^{1})/\mathrm{Im}(\delta^{0})$.
\end{Thm}

We have the following result which is similar to \cite[Proposition 3.3]{LRW}.

\begin{Prop}\label{parallel-paths-in-monomial-algebras}
Let $A=kQ_A/I_A$ be a monomial algebra and let $B=kQ_B/I_B$ be a (monomial) subalgebra of $A$ obtained by gluing two arrows $\alpha: e_1 \to e_2$ and $\beta: e_{n-1} \to e_n$ of $Q_A$. Let $\varphi:Q_{A} \to Q_{B}$ be the quiver morphism defined as in Remark \ref{quiver-morphism}. Then we have the following.

$(1)$ The map $\varphi: Q_{A} \to Q_{B}$ induces a surjective map $\widetilde{\varphi}: \mathcal{B}_A\rightarrow \mathcal{B}_B$ (denote $\widetilde{\varphi}(p)$ by $p^*$ for $p\in \mathcal{B}_A$) such that $\widetilde{\varphi}^{-1}(a^*)=\{a\}$ for $a^*\in (Q_B)_1\backslash \{\gamma^* \}$, $\widetilde{\varphi}^{-1}(\gamma^*)=\{\alpha,\beta\}$, $\widetilde{\varphi}^{-1}(f_i)=\{e_i\}$ for $3 \le i \le n-2$, $\widetilde{\varphi}^{-1}(f_1)=\{e_1,e_{n-1}\}$ and $\widetilde{\varphi}^{-1}(f_2)=\{e_2,e_{n}\}$.

$(2)$ Let $p,q \in \mathcal{B}_A$. Then $p \| q$ in $Q_A$ implies $p^{*} \| q^{*}$ in $Q_B$.

$(3)$ The map $\widetilde{\varphi}: \mathcal{B}_A\rightarrow \mathcal{B}_B$ induces $k$-linear maps $\psi_0: k((Q_A)_0\| \mathcal{B}_A) \to k((Q_B)_0\| \mathcal{B}_B)$, $\psi_1: k((Q_A)_1\| \mathcal{B}_A) \to k((Q_B)_1\| \mathcal{B}_B)$, $\psi_2: k(Z_A\| \mathcal{B}_A) \to k(Z_B\| \mathcal{B}_B)$.
\end{Prop}

Note that the definitions of the above maps $\psi_i$ ($0\leq i\leq 2$) are obvious, for example, the map $\psi_1: k((Q_A)_1 \| \mathcal{B}_A)\to k((Q_B)_1 \| \mathcal{B}_B)$ assigns $a \| p$ to $a^* \| p^*$ for each $a \| p \in k((Q_A)_1 \| \mathcal{B}_A)$.

From now on, we fix $A=kQ_A/I_A$ and $B=kQ_B/I_B$ to be the monomial algebras as in Proposition \ref{parallel-paths-in-monomial-algebras}, where $I_A=\mathrm{Span}(Z_A)$ and $I_B=\mathrm{Span}(Z_B)$. Then, by obvious identification, we have $Z_B=Z_A\cup Z_{new}$ and the following diagram:

\begin{center}
		\begin{tikzcd}
		0 \arrow[r] & k((Q_A)_{0}\| \mathcal{B}_A) \arrow[d, "\psi_{0}"] \arrow[r, "\delta_A^{0}"] & k((Q_A)_{1}\| \mathcal{B}_A) \arrow[d, "\psi_{1}"] \arrow[r, "\delta_A^{1}"] & k(Z_A\| \mathcal{B}_A) \arrow[d, "\psi_{2}"]\\
        0 \arrow[r] & k((Q_B)_{0}\| \mathcal{B}_B) \arrow[r, "\delta_B^{0}"] & k((Q_B)_{1}\| \mathcal{B}_B) \arrow[r, "\delta_B^{1}"] & k(Z_B\| \mathcal{B}_B)~
		\end{tikzcd}	
	    \quad \quad $(*)$
\end{center}
where $\psi_0: k((Q_A)_0\| \mathcal{B}_A) \to k((Q_B)_0\| \mathcal{B}_B)$, $\psi_1: k((Q_A)_1\| \mathcal{B}_A) \to k((Q_B)_1\| \mathcal{B}_B)$, $\psi_2: k(Z_A\| \mathcal{B}_A) \to k(Z_B\| \mathcal{B}_B)$ are the induced $k$-linear maps from the quiver morphism $\varphi: Q_{A} \to Q_{B}$ as mentioned in Proposition \ref{parallel-paths-in-monomial-algebras}. Note that the top and the bottom  complexes are truncations of the complexes $\mathcal{C}_{mon}$ of $A$ and of $B$, respectively.
Both squares in the diagram $(*)$ are not commutative in general, however, there are  close connections between the coboundary elements (resp. the cocycle elements) of the top complex and the coboundaries (resp. the cocycles) of the bottom complex in the diagram $(*)$.

In the next three sections, we will study the behaviour of Hochschild cohomology of degree one and zero and of the $\pi_1$-rank in case of gluing  two arbitrary arrows from a finite dimensional monomial algebra. Analogous to what we do in gluing two idempotents, we will consider two types of gluing: when we glue from the same block or from two different blocks.

\section{First Hochschild cohomology}\label{firstHoch}
In this section, we study the relation between the Lie algebras $\mathrm{HH}^1(A)$ and $\mathrm{HH}^1(B)$ when we glue two arrows $\alpha:e_1 \to e_2$ and $\beta:e_{n-1} \to e_{n}$ of $Q_A$, where $A$ is a monomial algebra. Recall that according to our assumption, the four vertices $e_1, e_2, e_{n-1}, e_{n}$ are pairwise different.

Firstly, in order to compare $\mathrm{Im}(\delta_A^0)$ with $\mathrm{Im}(\delta_B^0)$ we introduce the notion of special path with respect to gluing two arrows, which is a variation of special path with respect to gluing two idempotents.

\begin{Def} \label{Special-path}  {\rm(cf. \cite[Definition 3.6]{LRW})}
	Let $A=kQ_A/I_A$ be a monomial algebra and let $B=kQ_B/I_B$ be the algebra obtained from $A$ by gluing two arrows $\alpha$ and $\beta$. Let $p$ (resp. $q$) be a path in $\mathcal{B}_A$ of length $\ge 1$ either from $e_1$ to $e_{n-1}$ or from $e_{n-1}$ to $e_1$ (resp. either from $e_2$ to $e_n$ or from $e_n$ to $e_2$). We call $p$ (resp. $q$) a special path between $e_1$ and $e_{n-1}$ (resp. $e_2$ and $e_n$) if $\delta_B^0(f_1 \| p^*) \neq 0$ (resp. $\delta_B^0(f_2 \| q^*) \neq 0$), or equivalently, if there exists some $a\in (Q_A)_1$ such that $pa \notin I_A$ or $ap \notin I_A$(resp. $aq \notin I_A$ or $qa \notin I_A$).

\end{Def}

\begin{Notation} We denote by:
	\begin{itemize}
	\item $\mathrm{Sp}(1,n-1)$ (resp. $\mathrm{Sp}(2,n)$) the set of special paths between $e_1$ and $e_{n-1}$ (resp. $e_2$ and $e_n$) in $Q_A$, 
	\item  $\mathrm{Sp}(\alpha,\beta):=\mathrm{Sp}(1,n-1) \cup \mathrm{Sp}(2,n)$. We call an element in $\mathrm{Sp}(\alpha,\beta)$ a special path between  $\alpha$ and $\beta$, 
	\item $Z_{sp}(\alpha,\beta)$ the $k$-subspace of $\mathrm{Im}(\delta^0_B)$ generated by the elements $\delta_B^0(f_1 \| p^*)$ and $\delta_B^0(f_2 \| q^*)$, where $p \in \mathrm{Sp}(1,n-1)$ and $q \in \mathrm{Sp}(2,n)$,
	\item $\mathrm{sp}(\alpha,\beta)$ the dimension of the k-vector space $Z_{sp}(\alpha,\beta),$ 
	\item $\mathrm{sp}(1,n-1)$ (resp. $\mathrm{sp}(2,n)$) the cardinality of the set $\Sp(1,n-1)$ (resp. $\Sp(2,n)$).
	\end{itemize}
\end{Notation}
Note that this definition does not depend on whether the algebra $A$ is indecomposable or not. It is clear that if $A$ is a radical square zero algebra, then $\mathrm{sp(\alpha,\beta)}=0$. Moreover, when gluing a source arrow and a sink arrow, this definition can be simplified using the following notion:

\begin{Def} \label{Crucial-path}
Let $A=kQ_A/I_A$ be a monomial algebra and let $B=kQ_B/I_B$ be the algebra obtained from $A$ by gluing a source arrow $\alpha:e_1 \to e_2$ and a sink arrow $\beta:e_{n-1} \to e_n$. Let $\widetilde{p}$ be a path from $e_2$ to $e_{n-1}$ in $\mathcal{B}_A$ of length $\ge 1$. If $\beta\widetilde{p}\alpha \notin I_A$, then we say that $\widetilde{p}$ is a crucial path from $e_2$ to $e_{n-1}$. Denote by $\mathrm{Cp}(2,n-1)$ the set of all crucial paths and by $\mathrm{cp}(2,n-1)$ the cardinality of the set $\mathrm{Cp}(2,n-1)$.
\end{Def}

It is clear that if $\widetilde{p}$ is a crucial path from $e_2$ to $e_{n-1}$, then $p=\widetilde{p}\alpha$ is a special path from $e_1$ to $e_{n-1}$ and $q=\beta\widetilde{p}$ is a special path from $e_2$ to $e_n$. This implies that in gluing source-sink situation, we have
	\begin{equation*}
		\begin{split}
		\mathrm{Sp}(\alpha,\beta) & = \{\widetilde{p}\alpha,\beta\widetilde{p} \ |\  \widetilde{p} \in \mathrm{Cp}(2,n-1)\},\\
Z_{sp}(\alpha,\beta) & = \langle \delta_B^0(f_1 \| \widetilde{p}^*\gamma^*),\delta_B^0(f_2 \| \gamma^*\widetilde{p}^*) \ |\ \widetilde{p} \in \mathrm{Cp}(2,n-1) \rangle \\ & =	\langle \gamma^* \| \gamma^*\widetilde{p}^*\gamma^* \ |\ \widetilde{p} \in \mathrm{Cp}(2,n-1) \rangle,
		\end{split}
	\end{equation*}
	where the last line holds since $\delta_B^0(f_1 \| \widetilde{p}^*\gamma^*)=\gamma^*\|\gamma^*\widetilde{p}^*\gamma^*=-\delta_B^0(f_2 \| \gamma^*\widetilde{p}^*)$.

Note that in general $\mathrm{sp}(\alpha,\beta)$ is not equal to the cardinality of the set $\mathrm{Sp}(\alpha,\beta)$ ($=\mathrm{sp}(1,n-1)+\mathrm{sp}(2,n)$), which is not the case when gluing idempotents (compare to \cite[Remark 3.8]{LRW}). More precisely, we have the following remark:

\begin{Rem}\label{remark-on-Crucial-path}
As discussed above, in gluing source-sink situation we have $$\mathrm{sp}(\alpha,\beta)=\mathrm{dim}\,Z_{sp}(\alpha,\beta)=\mathrm{cp}(2,n-1)=\mathrm{sp}(1,n-1)=\mathrm{sp}(2,n),$$
which is the half of the cardinality of the set $\mathrm{Sp}(\alpha,\beta)$.
\end{Rem}

It is clear that if the source arrow $\alpha$ and the sink arrow $\beta$ are from different blocks of $A$, then $\mathrm{sp(\alpha,\beta)}=0$ since $\mathrm{cp}(2,n-1)=0$. However, if the source arrow $\alpha$ and the sink arrow $\beta$ are from the same block of $A$ and if there is a crucial path from $e_2$ to $e_{n-1}$ in $Q_A$, then $\mathrm{sp}(\alpha,\beta)$ is not equal to zero. In fact, $\mathrm{sp(\alpha,\beta)}$ can be arbitrarily large since $\mathrm{sp(\alpha,\beta)}\geq \mathrm{sp}(1,n-1)$ for gluing arbitrary two arrows and $\mathrm{sp(1,n-1)}$ can be arbitrarily large.

 The following proposition is a similar version of \cite[Proposition 3.9]{LRW}. It should be noted that,  although when gluing arrows and gluing idempotents the formulas of $\dim\,\Im(\delta_A^0)-\dim\,\Im(\delta_B^0)$ are similar, their actual values could be very different. For example, we will see in the next proposition that if we glue a source arrow and a sink arrow from the same block of $A$, then $\dim\,\Im(\delta_A^0)-\dim\,\Im(\delta_B^0)=2-\mathrm{sp(\alpha,\beta)}$, but if we glue a source idempotent and a sink idempotent from the same block of $A$, then $\dim\,\Im(\delta_A^0)-\dim\,\Im(\delta_B^0)=1$. We first recall the formal definitions of the maps $\delta^0_{(A)_0}$ and $\delta^0_{(A)_{\geq 1}}$ from \cite{LRW}.

 \begin{Def}\cite[Definition 3.4]{LRW} \label{the-decomposition-of-delta0}
We denote by $\delta^0_{(A)_0}$ (resp. $\delta^0_{(A)_{\geq 1}}$) the map $\delta_A^0$ restricted to the subspace $k((Q_A)_0\|(Q_A)_0)$ (resp. $k((Q_A)_0\|(\mathcal{B}_A)_{\geq 1})$). In particular, $\mathrm{Im}(\delta_{(A)_0}^0)$ is the $k$-vector space generated by the elements $\delta_A^0(e_i \| e_i)$ for $1\le i\le n$, and $\mathrm{Im}(\delta_{(A)_{\geq 1}}^0)$ is the $k$-vector space generated by the elements $\delta_A^0(e_i \| p)$ for $1\le i\le n$  where $p$ are paths of length at least $1$.
 \end{Def}
 
\begin{Prop}\label{Image-delta-zero-of-gluing-arrows} {\rm(Compare with \cite[Proposition 3.9]{LRW})}
	Let $A=kQ_A/I_A$ be a monomial algebra and let $B=kQ_B/I_B$ be the algebra obtained from $A$ by gluing two arbitrary arrows $\alpha:e_1 \to e_2$ and $\beta:e_{n-1} \to e_n$. Then $$\mathrm{dim}\,\mathrm{Im}(\delta_A^0)=\mathrm{dim}\,\mathrm{Im}(\delta_B^0)+2+c_B-c_A-\mathrm{sp(\alpha,\beta)},$$
where $c_A$ and $c_B$ are the number of connected components of $Q_A$ and $Q_B$ respectively.
	In particular, if we glue $\alpha$ and $\beta$ from the same block of $A$, then
	 $$\mathrm{dim}\,\mathrm{Im}(\delta_A^0)=\mathrm{dim}\,\mathrm{Im}(\delta_B^0)+2-\mathrm{sp(\alpha,\beta)};$$
	if we glue $\alpha$ and $\beta$ from two different blocks of $A$, then $$\mathrm{dim}\,\mathrm{Im}(\delta_A^0)=\mathrm{dim}\,\mathrm{Im}(\delta_B^0)+1.$$
\end{Prop}

\begin{proof}
	 Firstly we recall some notations from Section \ref{prelim}: the vertices of $Q_{A}$ are given by $e_1,e_2,\cdots,e_{n-1},e_n$ and the vertices of $Q_{B}$ are given by $f_1,f_2,\cdots,f_{n-2}$, where $f_1$ is obtained by gluing $e_1$ and $e_{n-1}$ and $f_2$ is obtained by gluing $e_2$ and $e_n$, and the arrows $\alpha$ and $\beta$ of $Q_A$ are identified as the arrow $\gamma^*$ of $Q_B$.
	
	We begin by describing the basis elements in $\mathrm{Im}(\delta_{A}^{0})$ and in $\mathrm{Im}(\delta_{B}^{0})$. As in the gluing idempotents case, we have the natural decomposition $\mathrm{Im}(\delta_A^0)=\mathrm{Im}(\delta_{(A)_0}^0) \oplus \mathrm{Im}(\delta_{(A)_{\ge 1}}^0)$ and $\mathrm{Im}(\delta_B^0)=\mathrm{Im}(\delta_{(B)_0}^0) \oplus \mathrm{Im}(\delta_{(B)_{\ge 1}}^0)$ as $k$-vector spaces, so it suffices to compare the $k$-subspaces $\mathrm{Im}(\delta_{(A)_0}^0)$ with $\mathrm{Im}(\delta_{(B)_0}^0)$,  and $\mathrm{Im}(\delta_{(A)_{\ge 1}}^0)$ with $\mathrm{Im}(\delta_{(B)_{\ge 1}}^0)$, respectively.
	
	(a1) We compare $\mathrm{Im}(\delta_{(A)_{\ge 1}}^0)$ with $\mathrm{Im}(\delta_{(B)_{\ge 1}}^0)$. First we observe that $\mathrm{Im}(\delta_{(A)_{\ge 1}}^0)$ is generated by the element $\delta_A^0(e_i \| p)$, where $p$ is a cycle at $e_i$ for $1\le i \le n$; similarly, $\mathrm{Im}(\delta_{(B)_{\ge 1}}^0)$ is generated by the element $\delta_B^0(f_i \| p^*)$, where $p^*$ is a cycle at $f_i$ for $1\le i \le n-2$. More specifically, 
	
	\begin{equation*}
		\begin{split}
		  \mathrm{Im}(\delta_{(B)_{\ge 1}}^0) & = \langle \delta_B^0(f_i \| p^*)\ |\ p\ is\ a\ cycle\ at\ e_i, 1\le i\le n \rangle\\ & \oplus \langle \delta_B^0(f_1 \| p^*),\delta_B^0(f_2 \| q^*) \ |\ p \in \mathrm{Sp}(1,n-1),q \in \mathrm{Sp}(2,n)\rangle \\ & = \langle \delta_B^0(f_i \| p^*)\ |\ p\ is\ a\ cycle\ at\ e_i, 1\le i\le n \rangle \oplus Z_{sp}(\alpha,\beta).
		\end{split}
	\end{equation*}
	From a direct computation, we have that $\delta_B^0(f_i \| p^*)=\psi_1(\delta_A^0(e_i \| p))$ for $p$ is a cycle at $e_i$ and $1\le i\le n$ (here we identify $f_{n-1}$ with $f_1$ and identify $f_n$ with $f_2$), which shows that we have a $k$-vector space decomposition of $\mathrm{Im}(\delta_{(B)_{\ge 1}}^0)$, that is, $$\mathrm{Im}(\delta_{(B)_{\ge 1}}^0)=\psi_1(\mathrm{Im}(\delta_{(A)_{\ge 1}}^0)) \oplus Z_{sp}(\alpha,\beta).$$ Hence $$\mathrm{dim}\,\mathrm{Im}(\delta_{(A)_{\ge 1}}^0)=\mathrm{dim}\,\mathrm{Im}(\delta_{(B)_{\ge 1}}^0)-\mathrm{sp(\alpha,\beta)}.$$

	(a2) We compare $\mathrm{Im}(\delta_{(A)_0}^0)$ with $\mathrm{Im}(\delta_{(B)_0}^0)$. It follows from \cite[Lemma 3.5]{LRW} that we have $\mathrm{dim}\,\mathrm{Im}(\delta_{(A)_0}^0)=\mathrm{dim}\,\mathrm{Im}(\delta_{(B)_0}^0)+2+c_B-c_A$. By combining this with the dimension formula  in (a1), we have $$\mathrm{dim}\,\mathrm{Im}(\delta_A^0)=\mathrm{dim}\,\mathrm{Im}(\delta_B^0)+2+c_B-c_A-\mathrm{sp(\alpha,\beta)}.$$
The statement follows.
\end{proof}

We have the following structural results when gluing a source arrow and a sink arrow.

\begin{Prop}\label{Image-delta-zero-of-gluing-source-and-sink-arrows}
	Let $A=kQ_A/I_A$ be a monomial algebra and let $B=kQ_B/I_B$ be the algebra obtained from $A$ by gluing a source arrow $\alpha:e_1 \to e_2$ and a sink arrow $\beta:e_{n-1} \to e_n$. Then we have the following equality as $k$-vector spaces: $$\langle \mathrm{Im}(\delta_B^0),\gamma^* \| \gamma^* \rangle = \psi_1(\mathrm{Im}(\delta_A^0)) \oplus Z_{sp}(\alpha,\beta),$$ where $\psi_1(\mathrm{Im}(\delta_A^0)) \simeq \mathrm{Im}(\delta_A^0)/\langle \alpha \| \alpha-\beta \| \beta \rangle$.
	In particular, if we glue $\alpha$ and $\beta$ from the same block of $A$, then
	$$\mathrm{Im}(\delta_B^0) \oplus \langle \gamma^* \| \gamma^* \rangle = \psi_1(\mathrm{Im}(\delta_A^0)) \oplus Z_{sp}(\alpha,\beta)$$ as $k$-vector spaces;
	if we glue $\alpha$ and $\beta$ from two different blocks of $A$, then $$\mathrm{Im}(\delta_B^0)=\psi_1(\mathrm{Im}(\delta_A^0)).$$
\end{Prop}

\begin{proof}
We first compare  $\mathrm{Im}(\delta_{(A)_0}^0)$ with  $\mathrm{Im}(\delta_{(B)_0}^0)$. 
 By definition, $\mathrm{Im}(\delta_{(A)_0}^0)$  is generated by the element $\delta_A^0(e_i \| e_i)$ for $1\le i\le n$. Since $\alpha:e_1 \to e_2$ is a source arrow and $\beta:e_{n-1} \to e_n$ is a sink arrow, then $\delta_A^0(e_1 \| e_1)=\alpha \| \alpha$, $$\delta_A^0(e_{n-1} \| e_{n-1})= \beta \| \beta - \sum_{t(b)=e_{n-1},b \in (Q_A)_1}b \| b,$$ and $\delta_A^0(e_n \| e_n)= - \beta \| \beta$, $$\delta_A^0(e_2 \| e_2)= \sum_{s(a)=e_2,a \in (Q_A)_1}a \| a - \alpha \| \alpha.$$ In addition, $$\delta_A^0(e_i \| e_i)= \sum_{s(a)=e_i,a \in (Q_A)_1} a \| a-\sum_{t(b)=e_i,b \in (Q_A)_1} b \| b$$ for $3\le i \le n-2$.
	
	Similarly, $\mathrm{Im}(\delta_{(B)_0}^0)$ is generated by the element of the form $\delta_B^0(f_i \| f_i)$ for $1\le i\le n-2$. By a direct computation, we have $$\delta_B^0(f_1 \| f_1)=\gamma^* \| \gamma^*-\sum_{t(b)=e_{n-1},b \in (Q_A)_1}b^* \| b^*=\psi_1(\delta_A^0(e_{n-1} \| e_{n-1})),$$ $$\delta_B^0(f_2 \| f_2)=\sum_{s(a)=e_2,a \in (Q_A)_1}a^* \| a^*-\gamma^* \| \gamma^*=\psi_1(\delta_A^0(e_2 \| e_2))$$ and $$\delta_B^0(f_i \| f_i)= \sum_{s(a)=e_i,a \in (Q_A)_1} a^* \| a^*-\sum_{t(b)=e_i,b \in (Q_A)_1} b^* \| b^*=\psi_1(\delta_A^0(e_i \| e_i))$$ for $3\le i \le n-2$.
	
	Therefore, we have $$\langle \mathrm{Im}(\delta_{(B)_0}^0),\gamma^* \| \gamma^* \rangle=\psi_1(\mathrm{Im}(\delta_{(A)_0}^0)).$$
	
	By combining the previous equality with the decomposition $\mathrm{Im}(\delta_{(B)_{\ge 1}}^0)=\psi_1(\mathrm{Im}(\delta_{(A)_{\ge 1}}^0)) \oplus Z_{sp}(\alpha,\beta)$ (cf. the proof of Proposition \ref{Image-delta-zero-of-gluing-arrows}) and noticing that $\psi_1$ in Diagram $(*)$ induces a natural map from $\mathrm{Im}(\delta_A^0)$ to $\langle \mathrm{Im}(\delta_B^0),\gamma^* \| \gamma^* \rangle$ with $\mathrm{Ker}(\psi_1)=\langle \alpha \| \alpha -\beta \| \beta \rangle$, we get $$\langle \mathrm{Im}(\delta_B^0),\gamma^* \| \gamma^* \rangle = \psi_1(\mathrm{Im}(\delta_A^0)) \oplus Z_{sp}(\alpha,\beta)$$ as $k$-vector spaces.
		
	Also note that, if we glue a source arrow $\alpha$ and a sink arrow $\beta$ from the same block of $A$, then by Lemma \ref{gammapi},
	$\gamma^* \| \gamma^* \notin \mathrm{Im}(\delta_B^0)$; if we glue a source arrow $\alpha$ and a sink arrow $\beta$ from two different blocks of $A$, then
	$\gamma^* \| \gamma^* \in \mathrm{Im}(\delta_B^0)$, $Z_{sp}(\alpha,\beta)=0$ and $c_B=c_A-1$. For the latter case, without loss of generality, we assume that $A$ has only two blocks, say $A_1$ and $A_2$. In addition, suppose that $\alpha:e_1 \to e_2$ is a source arrow in $A_1$ and $\beta:e_{n-1} \to e_{n}$ is a sink arrow in $A_2$. Let $I= \{1,2,\dots,n \}$, and denote the index set of the set of idempotents in $A_j$ by $I_j$ for $j=1,2$. Then $I$ is the disjoint union of $I_1$ and $I_2$. Note that (cf. \cite[Lemma 3.5]{LRW})
	\begin{equation*}
	\begin{split}
	0 & =\sum_{i \in I_2} \psi_1(\delta_A^0(e_i \| e_i)) \\ & =\psi_1(\delta_A^0(e_{n-1} \| e_{n-1}))+\psi_1(\delta_A^0(e_n \| e_n))+\sum_{i \in I_2 \backslash \{n-1,n\}} \psi_1(\delta_A^0(e_i \| e_i)) \\ & =\delta_B^0(f_1 \| f_1)-\gamma^* \| \gamma^*+\sum_{i \in I_2 \backslash \{n-1,n\}} \delta_B^0(f_i \| f_i),
	\end{split}
	\end{equation*}
	hence $\gamma^* \| \gamma^*=\delta_B^0(f_1 \| f_1)+\sum_{i \in I_2 \backslash \{n-1,n\}} \delta_B^0(f_i \| f_i)$ belongs to $\mathrm{Im}(\delta_{(B)_0}^0)$. We are done.
\end{proof}

\begin{Rem}\label{dim-Image}
\begin{itemize}	
\item[$(1)$]  When we glue two arbitrary arrows, the structural results in Proposition \ref{Image-delta-zero-of-gluing-source-and-sink-arrows} are not true in general. Specifically, $\alpha\|\alpha-\beta\|\beta$ may be not in $\Im(\delta_A^0)$ when we glue from the same block (cf. Example \ref{eg3}) or from different block (cf. Example \ref{eg5}). Moreover, $\Im(\delta_B^0)$ may be not a subset of $\psi_1(\Im(\delta_A^0))\oplus Z_{sp}(\alpha,\beta)$ when we glue from the same block (cf. Example \ref{eg3}), and $\Im(\delta_B^0)$ may be not equal to $\psi_1(\Im(\delta_A^0))$ when we glue from different blocks (cf. Example \ref{eg5}).
	
\item[$(2)$] Similar to the gluing idempotents situation, even if we glue a source arrow $\alpha$ and a sink arrow $\beta$ from the same block of $A$,  $\psi_1$ does not induce a map from $\Im(\delta_A^0)$ to $\Im(\delta_B^0)$ since $\psi_1(\delta_A^0(e_1\|e_1))=\gamma^* \| \gamma^* \notin \Im(\delta_B^0)$ (cf. Example \ref{eg1}). Nevertheless, the natural map $$\psi_1:\mathrm{Im}(\delta_A^0) \to \langle \mathrm{Im}(\delta_B^0),\gamma^* \| \gamma^* \rangle=\psi_1(\mathrm{Im}(\delta_A^0)) \oplus Z_{sp}(\alpha,\beta)$$ induces the following exact sequence of $k$-vector spaces:
	 $$0\to\langle \alpha \| \alpha - \beta \| \beta \rangle \to \Im(\delta_A^0) \stackrel{\psi_1}{\longrightarrow} \psi_1(Im(\delta_A^0)) \oplus Z_{sp}(\alpha,\beta)\to  Z_{sp}(\alpha,\beta)\to 0.$$
\end{itemize}
\end{Rem}

\begin{Cor}\label{Imange-glue-arrows-radical-square-zero}
Let $A=kQ_A/I_A$ be a radical square zero algebra and let $B=kQ_B/I_B$ be the algebra obtained from $A$ by gluing a source arrow $\alpha:e_1 \to e_2$ and a sink arrow $\beta:e_{n-1} \to e_n$. Then $$\langle \mathrm{Im}(\delta_B^0),\gamma^* \| \gamma^* \rangle = \psi_1(\mathrm{Im}(\delta_A^0))\simeq \mathrm{Im}(\delta_A^0)/\langle \alpha \| \alpha-\beta \| \beta \rangle$$ as $k$-vector spaces, and $$\mathrm{dim}\,\mathrm{Im}(\delta_A^0)=\mathrm{dim}\,\mathrm{Im}(\delta_B^0)+2+c_B-c_A.$$
	In particular, if we glue $\alpha$ and $\beta$ from the same block of $A$, then
	$\mathrm{Im}(\delta_B^0) \oplus \langle \gamma^* \| \gamma^* \rangle = \psi_1(\mathrm{Im}(\delta_A^0))$ as $k$-vector spaces, and $\mathrm{dim}\,\mathrm{Im}(\delta_A^0)=\mathrm{dim}\,\mathrm{Im}(\delta_B^0)+2;$
	if we glue $\alpha$ and $\beta$ from two different blocks of $A$, then $\mathrm{Im}(\delta_B^0)=\psi_1(\mathrm{Im}(\delta_A^0))$ and $\mathrm{dim}\,\mathrm{Im}(\delta_A^0)=\mathrm{dim}\,\mathrm{Im}(\delta_B^0)+1.$
\end{Cor}

\begin{proof}
In radical square zero case there is no special path between $e_1$ and $e_{n-1}$ (resp. between $e_2$ and $e_n$), hence $Z_{sp}(\alpha,\beta)=0$. The statement follows directly from Proposition \ref{Image-delta-zero-of-gluing-source-and-sink-arrows}.
\end{proof}

We now proceed to compare $\mathrm{Ker}(\delta^1_A)$ with $\mathrm{Ker}(\delta^1_B)$. We first prove a similar version of \cite[Proposition 3.12]{LRW} when we glue a source arrow and a sink arrow.

\begin{Prop}\label{Kernel-delta-one}
	Let $A=kQ_A/I_A$ be a (not necessarily indecomposable) monomial algebra and let $B=kQ_B/I_B$ be the algebra obtained from $A$ by gluing a source arrow $\alpha:e_1 \to e_2$ and a sink arrow $\beta:e_{n-1} \to e_n$. Then there exists a (restricted) Lie algebra homomorphism $\mathrm{Ker}(\delta^1_A) \to \mathrm{Ker}(\delta^1_B)$ induced from $\psi_1: k((Q_A)_1\| \mathcal{B}_A) \to k((Q_B)_1\| \mathcal{B}_B)$, with kernel generated by the element $\alpha \| \alpha - \beta \| \beta$, which we still denote by $\psi_1$.
\end{Prop}

\begin{proof}
	First notice that when gluing source arrow $\alpha:e_1 \to e_2$ and sink arrow $\beta:e_{n-1} \to e_n$, we have that $Z_{new}=\{\xi^*\gamma^*\eta^*\ |\ t(\eta)=e_{n-1},s(\xi)=e_2;\eta,\xi \in (Q_A)_1\}$. Hence for $I_A=\langle Z_A \rangle$, by the obvious identification we can write $Z_B=Z_A \cup Z_{new}$ such that $I_B=\langle Z_B \rangle$.
	
	Let $a \| p \in k((Q_A)_{1}\| \mathcal{B}_A)$ and assume that $p=a_m \dots a_1$. Let $\psi_1(a \| p)=a^{*} \| p^{*}$ be the corresponding element in $k((Q_B)_1\| \mathcal{B}_B)$. On the one hand, we have
	$$\psi_2(\delta_{A}^{1}(a \| p))=\psi_2(\sum_{r \in Z_A} r\| r^{a \| p})=\sum_{r \in Z_A} r\| r^{a^* \| p^*};$$
	On the other hand, we have
	$$\delta_{B}^{1}(\psi_1(a \| p))=\delta_{B}^{1}(a^* \| p^*)=\sum_{r \in Z_A} r\| r^{a^* \| p^*}+\sum_{r' \in Z_{new}} r'\| {r'}^{a^* \| p^*}.$$
	
	We consider two cases.
	
	(b1) If $a$ is not a loop arrow and $p$ is a path parallel to $a$. In this case it is easy to see that once $a^*$ appears in some $r'\in Z_{new}$ (note that $r'\notin \mathcal{B}_B$), the element obtained from $r'$ by replacing $a^*$ by $p^*$ is again not in $\mathcal{B}_B$, whence have $\sum_{r' \in Z_{new}} r'\| {r'}^{a^* \| p^*}=0$. In fact, since $r'$ has the form $\xi^*\gamma^*\eta^*$, where $\eta,\xi$ are arrows in $(Q_A)_1$ with $t(\eta)=e_{n-1}$ and $s(\xi)=e_2$, then $r'^{a^* \| p^*} \ne 0$ only if $a^*=\eta^*,\gamma^*$ or $\xi^*$. If $a^*=\eta^*$,  then since $a \| p$ we get $t(a_m)=t(p)=t(a)=t(\eta)=e_{n-1}$. This implies that $r'^{a^* \| p^*}=\xi^*\gamma^*a_m^* \dots a_1^*=0$ in $B$ since $\xi^*\gamma^*a_m^* \in Z_{new}$. The left two cases are similar. Therefore $\psi_2(\delta_{A}^{1}(a \| p))=\delta_{B}^{1}(\psi_1(a \| p))$.
	
	(b2) If $a$ is a loop arrow at $e_i$ and $p=e_i$ or $p$ is an oriented cycle at $e_i$ for $1 \le i \le n$. Actually, $3 \le i \le n-2$ since $\alpha:e_1 \to e_{n-1}$ is a source arrow and $\beta:e_2 \to e_n$ is a sink arrow. Then $\sum_{r' \in Z_{new}} r'\| {r'}^{a^* \| p^*}=0$ in $B$ since $a^*$ does not appear in any $r'\in Z_{new}$. Therefore $\psi_2(\delta_{A}^{1}(a \| p))=\delta_{B}^{1}(\psi_1(a \| p))$.
	
	As a result, we get a $k$-linear map $\psi_{1}: \mathrm{Ker}(\delta_{A}^{1}) \to \mathrm{Ker}(\delta_{B}^{1})$ induced from the following mapping: $a \| e_i \mapsto a^* \| f_i$ ($i\neq 1,2,n-1,n$), $\alpha\|\alpha \mapsto \gamma^*\| \gamma^*$, $\beta\|\beta \mapsto \gamma^*\|\gamma^*$, $a \| p\mapsto a^* \| p^*$ ($p \in \mathcal{B}_A$ has length $\geq 1$). It is also clear that $\psi_{1}: \mathrm{Ker}(\delta_{A}^{1}) \to  \mathrm{Ker}(\delta_{B}^{1})$ has kernel generated by $\alpha \| \alpha - \beta \| \beta$ and preserves the Lie bracket, since $\psi_1$ preserves the parallel paths.
\end{proof}

In order to prove a general version of Proposition \ref{Kernel-delta-one} for gluing two arbitrary arrows, we need a similar characteristic condition as in \cite[Assumption 3.11]{LRW}.

\begin{Assumption}\label{assum}
Let $A=kQ_A/I_A$ be a monomial algebra and let $B=kQ_B/I_B$ be the algebra obtained from $A$ by gluing two arrows $\alpha:e_1 \to e_2$ and $\beta:e_{n-1} \to e_n$. For each loop $a$ at $e_1,e_2,e_{n-1}$ or $e_n$ with $a^m\in Z_A$ for some $m\geq 2$, we have that $\mathrm{char}(k)\nmid m$.
\end{Assumption}

Clearly, Assumption \ref{assum} holds when the characteristic of the field $k$ is zero or big enough.

\begin{Prop}\label{Kernel-delta} {\rm(Compare to \cite[Proposition 3.12]{LRW})}
	Let $A=kQ_A/I_A$ be a monomial algebra and let $B=kQ_B/I_B$ be the algebra obtained from $A$ by  gluing  two arbitrary arrows  $\alpha:e_1 \to e_2$ and $\beta:e_{n-1} \to e_n$. Then under Assumption \ref{assum}, there exists a (restricted) Lie algebra homomorphism $\psi_1:\mathrm{Ker}(\delta^1_A) \to \mathrm{Ker}(\delta^1_B)$ induced from $\psi_1: k((Q_A)_1\| \mathcal{B}_A) \to k((Q_B)_1\| \mathcal{B}_B)$ with kernel generated by the element $\alpha \| \alpha - \beta \| \beta$.
\end{Prop}

\begin{proof}
The proof is similar to Proposition \ref{Kernel-delta-one} but considers as part of case (b2) also the case when $a$ is a loop at $e_i$ for $i=1,2,n-1$ or $n$. If $i=1$, that is, $a$ is a loop at $e_1$, then there is a relation $r=a^s$ in $Z_A$ for some $s\geq 2$ since $A$ is finite dimensional. It follows that $sr\|a^{s-1}$ is a summand of $\delta_A^1(a\|e_1)$ which cannot be cancelled by other summands. Whence $\delta_A^1(a\|e_1)=0$ implies that $\mathrm{char}(k)|s$. Therefore if $\mathrm{chak}(k) \nmid s$, then $a\|e_1\notin \mathrm{Ker}(\delta_A^1)$ and $\psi_1(a\|e_1) \notin \mathrm{Ker}(\delta_B^1)$. The left cases are similar, and under Assumption \ref{assum} we can exclude that $a\|e_i$ belongs to $\mathrm{Ker}(\delta_A^1)$ for $i=1,2,n-1$ and $n$. The statement follows.
\end{proof}

Note that the map $\psi_1$ from $\mathrm{Ker}(\delta^1_A)$ to $\mathrm{Ker}(\delta^1_B)$ is injective when we glue two idempotents of $A$, but in gluing two arrows case $\psi_1$ changes to a map with one-dimensional kernel.

In the following Remark, the first item is a similar version of \cite[Remark 3.13]{LRW}, and the second item is a variation of \cite[Remark 3.14 (1)]{LRW}.

\begin{Rem}
	\begin{itemize}
		\item [(i)]  If there is no loop at $e_1,e_2,e_{n-1}$ or $e_n$, then we do not need Assumption \ref{assum} when we glue  two arbitrary arrows $\alpha:e_1\to e_2$ and $\beta:e_{n-1}\to e_n$. Indeed, the characteristic condition only makes sense when there are loops at $e_1,e_2,e_{n-1}$ or $e_n$.
		\item [(ii)] If $A$ (hence also $B$) is radical square zero, then we do not need Assumption \ref{assum}. Indeed, in this case Assumption \ref{assum} is equivalent to $\mathrm{char}(k)\neq 2$. However, the loop $a$ at  $e_1,e_2,e_{n-1}$ or $e_n$ must appear in a relation of the form $\alpha a, a\alpha,\beta a$ or $a\beta$, which gives rise to $a\|e_i\notin \mathrm{Ker}(\delta_A^1)$ and $\psi_1(a\|e_1) \notin \mathrm{Ker}(\delta_B^1)$ for $i=1,2,n-1$ or $n$.
	\end{itemize}
\end{Rem}

In order to describe the elements in $\mathrm{Ker}(\delta_{B}^{1})$ which are in the complement of the subspace $\psi_1(\mathrm{Ker}(\delta_{A}^{1}))$, we introduce some further notation.

\begin{Def}\label{Special-pair} \rm(Compare to \cite[Definition 3.15]{LRW})
	Let $A=kQ_A/I_A$ be a (not necessarily indecomposable) monomial algebra and let $B=kQ_B/I_B$ be the algebra obtained from $A$ by gluing the arrow $\alpha:e_1 \to e_2$ and the arrow $\beta:e_{n-1} \to e_n$. Let $a$ be an arrow and $p$ be a path in $\mathcal{B}_A$. We say that $(a,p)$ is a special pair with respect to the gluing of $\alpha$ and $\beta$ if the following four conditions are satisfied:
	
	(1) the starting vertex or the ending vertex of $a$ is $e_i$ for $i=1,2,n-1$ or $n$;
	
	(2) $a \nparallel p$ holds in $Q_A$;
	
	(3) $a^* \| p^*$ holds in $Q_B$ but $a^* \| p^*$ is not equal to $\gamma^* \| \gamma^*$;
	
	(4) if $a=\alpha$ (resp. $a=\beta$), then $p\nparallel \beta$ (resp. $p\nparallel \alpha$); if $p=\alpha$ (resp. $p=\beta$), then $a \nparallel \beta$ (resp. $a \nparallel \alpha$).
\end{Def}

	\begin{Notation} We denote by:
	\begin{itemize}
	\item $\mathrm{Spp}(\alpha,\beta)$ the set of all special pairs with respect to the gluing of $\alpha$ and $\beta$,
	\item   $\langle\mathrm{Spp}(\alpha,\beta)\rangle$ the $k$-subspace of $k((Q_B)_1\| \mathcal{B}_B)$ generated by the elements $a^* \| p^*$, where $(a,p) \in \mathrm{Spp}(\alpha,\beta)$, 
	\item  $Z_{spp}(\alpha,\beta)$ the intersection of $\langle\mathrm{Spp}(\alpha,\beta)\rangle$ and $\mathrm{Ker}(\delta^1_B)$, 
	\item $\mathrm{kspp}(\alpha,\beta)$ the dimension of the $k$-subspace $Z_{spp}(\alpha,\beta)$ of $\mathrm{Ker}(\delta^1_B)$. 
	\end{itemize}
	\end{Notation}

\begin{Rem}
\begin{itemize}
		\item [(i)] The reason why we just consider the case for $i=1,2,n-1$ or $n$ in Definition \ref{Special-pair} follows from the fact that if the $s(a)$ or $t(a)$ equals to $e_i$ for $3 \le i \le n-2$, then $a^* \| p^*$ implies $a \| p$ which shows that, if $a^* \| p^* \in \mathrm{Ker}(\delta_B^1)$, then $a^* \| p^*$ is not in the complement of $\psi_1(\mathrm{Ker}(\delta_{A}^{1}))$. In Condition (3) of Definition \ref{Special-pair} we require that $a^* \| p^* \neq \gamma^* \| \gamma^*$ since $\gamma^* \| \gamma^*=\psi_1(\alpha \| \alpha)$ is an element in $\psi_1(\mathrm{Ker}(\delta_{A}^{1}))$. Note that \cite[Definition 3.15]{LRW} has three conditions which correspond to the first three conditions in Definition \ref{Special-pair}. Example \ref{eg3} explains why we request Condition (4) in Definition \ref{Special-pair}.

\item [(ii)] The cardinality of the set $\mathrm{Spp}(\alpha,\beta)$ can be arbitrarily large even for radical square zero algebras (cf. Example \ref{eg6}), although in this case $\mathrm{Sp}(\alpha,\beta)=\emptyset$.
\end{itemize}
\end{Rem}

\begin{Rem}\label{special-pair-gluing-source-sink}
\begin{itemize}
		\item [(i)] Note that Definition \ref{Special-pair} holds in the case of gluing  two arbitrary arrows and can be simplified when gluing  a source arrow $\alpha$  and  a sink arrow $\beta$. In fact, if we glue a source arrow and a sink arrow, then the conditions (1)-(3) imply the condition (4). Moreover, it is not difficult to see that in this case $\mathrm{Spp}(\alpha,\beta)=\{(\alpha,\beta p \alpha),(\beta,\beta p \alpha) \ |\ p \in \mathrm{Cp}(2,n-1)\}$, where $\mathrm{Cp}(2,n-1)$ is defined in Definition \ref{Crucial-path}. Hence $Z_{spp}(\alpha,\beta)=0$ when the source arrow $\alpha$ and the sink arrow $\beta$ are from two different blocks since in this case $\mathrm{Cp}(2,n-1)=\emptyset$.
	
\item [(ii)] The above observation implies that when gluing a source arrow $\alpha:e_1 \to e_2$ and a sink arrow $\beta:e_{n-1} \to e_n$ we have that $\langle\mathrm{Spp}(\alpha,\beta)\rangle=\langle \gamma^* \| \gamma^*p^*\gamma^* \ |\ p\in \mathrm{Cp}(2,n-1) \rangle=Z_{sp}(\alpha,\beta) \subset \mathrm{Im}(\delta_B^0) \subset \mathrm{Ker}(\delta_B^1)$. Hence $Z_{spp}(\alpha,\beta)=\langle\mathrm{Spp}(\alpha,\beta)\rangle \cap \mathrm{Ker}(\delta_B^1)=\langle\mathrm{Spp}(\alpha,\beta)\rangle=Z_{sp}(\alpha,\beta)$ (note that in general $Z_{sp}(\alpha,\beta) \subset Z_{spp}(\alpha,\beta)$) and $\mathrm{kspp(\alpha,\beta)=sp(\alpha,\beta)} =\mathrm{cp}(2,n-1)$. This shows that in this case $Z_{spp}(\alpha,\beta)$ is generated by the elements of the form $a^* \| p^*$. However in general, a generator of $Z_{spp}(\alpha,\beta)$ is just a $k$-linear combination of the elements of the form $a^* \| p^*$ (cf. Example \ref{eg2}).
\end{itemize}
\end{Rem}

\begin{Prop}\label{Kernel-delta-one-of-gluing-source-sink-arrows}
	Let $A=kQ_A/I_A$ be a (not necessarily indecomposable) monomial algebra and let $B=kQ_B/I_B$ be the algebra obtained from $A$ by gluing a source arrow $\alpha:e_1 \to e_2$ and a sink arrow $\beta:e_{n-1} \to e_n$. Then $$\mathrm{Ker}(\delta_B^1)=\psi_1(\mathrm{Ker}(\delta_A^1)) \oplus Z_{spp}(\alpha,\beta)$$ as $k$-vector spaces, where $$\psi_1(\mathrm{Ker}(\delta_A^1)) \simeq \mathrm{Ker}(\delta_A^1)/\langle \alpha \| \alpha - \beta \| \beta \rangle \mbox{ and }Z_{spp}(\alpha,\beta)=Z_{sp}(\alpha,\beta)=\langle \gamma^* \| \gamma^*p^*\gamma^* \ |\ p\in \mathrm{Cp}(2,n-1) \rangle.$$ In particular, we have $\mathrm{dim}\,\mathrm{Ker}(\delta_B^1)=\mathrm{dim}\,\mathrm{Ker}(\delta_A^1)-1+\mathrm{cp}(2,n-1)$.
	
	 In particular, if we glue $\alpha$ and $\beta$ from the same block of $A$, then there is a decomposition as $k$-vector spaces, $$\mathrm{Ker}(\delta_B^1)=\psi_1(\mathrm{Ker}(\delta_A^1)) \oplus \langle \gamma^* \| \gamma^*p^*\gamma^* \ |\ p\in \mathrm{Cp}(2,n-1) \rangle, $$ and we have $\mathrm{dim}\,\mathrm{Ker}(\delta_B^1)=\mathrm{dim}\,\mathrm{Ker}(\delta_A^1)-1+\mathrm{cp}(2,n-1);$
	 if we glue $\alpha$ and $\beta$ from two different blocks of $A$, then $$\mathrm{Ker}(\delta_B^1)=\psi_1(\mathrm{Ker}(\delta_A^1))$$ and $\mathrm{dim}\,\mathrm{Ker}(\delta_B^1)=\mathrm{dim}\,\mathrm{Ker}(\delta_A^1)-1.$
\end{Prop}

\begin{proof}
	By Proposition $\ref{Kernel-delta-one}$, we only need to describe the elements $\theta$ in $\mathrm{Ker}(\delta_{B}^{1})$ which are in the complement of the subspace $\psi_1(\mathrm{Ker}(\delta_{A}^{1}))$. According to the proof of Proposition $\ref{Kernel-delta-one}$, we may assume that $\theta$ is a linear combination of the elements of the form $a^* \| p^*$ such that $(a,p)$ is a special pair with respect to the gluing of $\alpha$ and $\beta$. Clearly in this case $\theta\in Z_{spp}(\alpha,\beta)$, where $Z_{spp}(\alpha,\beta)$ is the subspace of $\mathrm{Ker}(\delta_{B}^{1})$ defined in Definition \ref{Special-pair}. Therefore, we have the following decomposition: $\mathrm{Ker}(\delta_B^1)=\psi_1(\mathrm{Ker}(\delta_A^1)) \oplus Z_{spp}(\alpha,\beta)$. Hence the dimension formula follows. The second statement follows from Remark $\ref{special-pair-gluing-source-sink}$.
\end{proof}

%\begin{Rem}\label{kernel-delta-one-of-gluing-arbitrary-arrows}
	We have a similar version in the case of gluing of two arbitrary  arrows under the characteristic condition Assumption \ref{assum}. More specifically,  all the results in Proposition $\ref{Kernel-delta-one-of-gluing-source-sink-arrows}$ hold when gluing two arbitrary arrows except the structure of $Z_{spp}(\alpha,\beta)$. In fact,  in general $Z_{spp}(\alpha,\beta)$ strictly contains $Z_{sp}(\alpha,\beta)$ (cf. Example \ref{eg4}).
%\end{Rem}

\begin{Prop}\label{Kernel-delta-one-structure} {\rm(Compare with \cite[Proposition 3.17]{LRW})}
	Let $A=kQ_A/I_A$ be a monomial algebra and let $B=kQ_B/I_B$ be the algebra obtained from $A$ by gluing two arrows $\alpha:e_1 \to e_2$ and $\beta:e_{n-1} \to e_n$. Then under Assumption \ref{assum}, $$\mathrm{Ker}(\delta_B^1)=\psi_1(\mathrm{Ker}(\delta_A^1)) \oplus Z_{spp}(\alpha,\beta)$$ as $k$-vector spaces, where $\psi_1(\mathrm{Ker}(\delta_A^1)) \simeq \mathrm{Ker}(\delta_A^1)/\langle \alpha \| \alpha - \beta \| \beta \rangle $. In particular, we have $\dim\Ker(\delta_B^1)=\dim\Ker(\delta_A^1)-1+\mathrm{kspp(\alpha,\beta)}$.
\end{Prop}

\begin{proof}
	Under Assumption \ref{assum}, a similar proof of Proposition \ref{Kernel-delta-one-of-gluing-source-sink-arrows} shows that the elements in $\mathrm{Ker}(\delta_{B}^{1})$ which are in the complement of the subspace $\psi_1(\mathrm{Ker}(\delta_{A}^{1}))$ are exactly in $Z_{spp}(\alpha,\beta)$. Therefore the proof follows from Proposition \ref{Kernel-delta}.
\end{proof}

\begin{Rem}\label{exact-seq-of-Ker}
	By Proposition $\ref{Kernel-delta-one-structure}$, and under Assumption \ref{assum}, if we glue  two arbitrary arrows of $A$, then there is a canonical exact sequence of $k$-vector spaces as follows:
	$$0\to \langle \alpha \| \alpha - \beta \| \beta \rangle\to \Ker(\delta_A^1) \stackrel{\psi_1}{\longrightarrow} \Ker(\delta_B^1)\to Z_{spp}(\alpha,\beta)\to 0.$$
In particular, when gluing a source arrow $\alpha$ and a sink arrow $\beta$, the characteristic condition is not necessary, and $Z_{spp}(\alpha,\beta)$ is exactly $Z_{sp}(\alpha,\beta)$ and vanishes when $\alpha$ and $\beta$ are from different blocks of $A$.
\end{Rem}

In the remaining part of this section, we discuss the relationship between the Lie algebras $\HH^1(A)$ and $\HH^1(B)$ under the operation of gluing two arrows. In the following, we just write the image of $\gamma^*\|\gamma^*$ in $\mathrm{HH}^{1}(B)$ as $\gamma^*\|\gamma^*$ for simplicity.

\begin{Thm}\label{hh1-gluing-source-sink-arrows} {\rm(Compare with \cite[Corollary 3.24]{LRW})}
	Let $A=kQ_A/I_A$ be a monomial algebra and let $B=kQ_B/I_B$ be the algebra obtained from $A$ by gluing a source arrow $\alpha:e_1 \to e_2$ and a sink arrow $\beta:e_{n-1} \to e_n$. Then $$\mathrm{HH}^{1}(A)\simeq \mathrm{HH}^{1}(B)/I$$ as Lie algebras, where $I:=Y/\mathrm{Im}(\delta_B^0)$ is a Lie ideal of $\mathrm{HH}^{1}(B)$ and $Y:=\psi_1(\mathrm{Im}(\delta_A^0)) \oplus Z_{sp}(\alpha,\beta)$. In particular, if we glue from the same block, then we have a Lie algebra isomorphism $$\mathrm{HH}^{1}(A)\simeq \mathrm{HH}^{1}(B)/\langle \gamma^*\|\gamma^* \rangle;$$
	if we glue from different blocks, then there is a Lie algebra isomorphism $$\mathrm{HH}^{1}(A)\simeq \mathrm{HH}^{1}(B).$$
\end{Thm}

\begin{proof}
	By Proposition $\ref{Kernel-delta-one}$, there exists a Lie algebra homomorphism $\psi_1: \mathrm{Ker}(\delta^1_A)\to \mathrm{Ker}(\delta^1_B)$, which is induced from the canonical map $\psi_1: k((Q_A)_1\| \mathcal{B}_A) \to k((Q_B)_1\| \mathcal{B}_B)$ ($a\| p\mapsto a^*\| p^*$). Define $Y:=\psi_1(\mathrm{Im}(\delta_A^0)) \oplus Z_{sp}(\alpha,\beta)$. Then Proposition \ref{Image-delta-zero-of-gluing-source-and-sink-arrows} implies that $Y=\mathrm{Im}(\delta_B^0)\oplus\langle \gamma^*\| \gamma^*\rangle$ when we glue from the same block, and $Y=\mathrm{Im}(\delta_B^0)$ when we glue from the different blocks. Consequently, we have that $I:=Y/\mathrm{Im}(\delta_B^0)$ is a one-dimensional Lie ideal of $\mathrm{HH}^{1}(B)$ generated by $\gamma^*\|\gamma^*$ when we glue from the same block, otherwise $I$ is zero.
	
	We claim that $Y$ is a Lie ideal of $\mathrm{Ker}(\delta_B^1)$ and $\mathrm{Im}(\delta_B^0)$ is a Lie ideal of $Y$. If we glue from the same block, then we show this in three steps, and if we glue from different blocks, then there is nothing to prove  since $Y=\mathrm{Im}(\delta_B^0)$.
	
	Firstly, we show that $Y$ is a Lie algebra. It suffices to show that $Y$ is closed under the Lie bracket. By the decomposition of $Y$ we have  $$[Y,Y]=[\mathrm{Im}(\delta_B^0),\mathrm{Im}(\delta_B^0)]+[\gamma^*\|\gamma^*,\mathrm{Im}(\delta_B^0)]+[\mathrm{Im}(\delta_B^0),\gamma^*\|\gamma^*]+[\gamma^*\|\gamma^*,\gamma^*\|\gamma^*],$$ it is clear that $[\gamma^*\|\gamma^*,\gamma^*\|\gamma^*]=0$ and $[\mathrm{Im}(\delta_B^0),\mathrm{Im}(\delta_B^0)] \subset \mathrm{Im}(\delta_B^0)$. We need only to consider $[\mathrm{Im}(\delta_B^0),\gamma^*\|\gamma^*]$. Note that $\gamma^*\|\gamma^* \in \mathrm{Ker}(\delta_B^1)$ and $\mathrm{Im}(\delta_B^0)$ is a Lie ideal of $\mathrm{Ker}(\delta_B^1)$, hence $[\mathrm{Im}(\delta_B^0),\gamma^*\|\gamma^*] \subset [\mathrm{Im}(\delta_B^0),\mathrm{Ker}(\delta_B^1)]$ $\subset \mathrm{Im}(\delta_B^0)$. Therefore, we have $[Y,Y] \subset \mathrm{Im}(\delta_B^0) \subset Y$, that is, $Y$ is a Lie algebra.
	
	Obviously, $\mathrm{Im}(\delta_B^0)$ is a Lie ideal of $Y$ since $[\mathrm{Im}(\delta_B^0),Y] \subset [\mathrm{Im}(\delta_B^0),\mathrm{Ker}(\delta_B^1)] \subset \mathrm{Im}(\delta_B^0)$.
	
	Finally, we show that $Y$ is a Lie ideal of $\mathrm{Ker}(\delta_B^1)$. Note that $[Y,\mathrm{Ker}(\delta_B^1)]=[\mathrm{Im}(\delta_B^0),\mathrm{Ker}(\delta_B^1)]+[\gamma^*\|\gamma^*,\mathrm{Ker}(\delta_B^1)]$ and $[\mathrm{Im}(\delta_B^0),\mathrm{Ker}(\delta_B^1)] \subset \mathrm{Im}(\delta_B^0) \subset Y$, hence it suffices to show that $[\gamma^*\|\gamma^*,\mathrm{Ker}(\delta_B^1)]$ is contained in $Y$. By Proposition \ref{Kernel-delta-one-of-gluing-source-sink-arrows}, we have $[\gamma^*\|\gamma^*,\mathrm{Ker}(\delta_B^1)]=[\gamma^*\|\gamma^*,\psi_1(\mathrm{Ker}(\delta_A^1))]+[\gamma^*\|\gamma^*,Z_{spp}(\alpha,\beta)]$ and $Z_{spp}(\alpha,\beta)=Z_{sp}(\alpha,\beta)=\langle \gamma^* \| \gamma^*p^*\gamma^* \ |\ p\in \mathrm{Cp}(2,n-1) \rangle $. It follows that $[\gamma^*\|\gamma^*,Z_{spp}(\alpha,\beta)] \subset Z_{sp}(\alpha,\beta) \subset Y$. We now claim that $[\gamma^*\|\gamma^*,\psi_1(\mathrm{Ker}(\delta_A^1))]=0$. In fact, it is enough to show that $[\gamma^*\|\gamma^*,a^*\|q^*]=a^*\|q^{*\ \gamma^*\|\gamma^*}-\delta_{a^*}^{\gamma^*}\gamma^*\|q^*$ is zero for any $a^*\|q^*=\psi_1(a\|q)$ with $a\|q$ appearing as a summand of some element in $\mathrm{Ker}(\delta_A^1)$, where $\delta_{a^*}^{\gamma^*}$ denotes the Kronecker symbol. If $a^*\|q^{*\ \gamma^*\|\gamma^*}=0$ but $[\gamma^*\|\gamma^*,a^*\|q^*]=-\delta_{a^*}^{\gamma^*}\gamma^*\|q^*\neq 0$, this yields that $a^*=\gamma^*$. It follows that $a=\alpha$ (resp. $\beta$), by combining $a\|q$ and $\alpha$ (resp. $\beta$) is a source (resp. sink) arrow, we have $a\|q=\alpha\|\alpha$ (resp. $a\|q=\beta\|\beta$). Now $a^*\|q^*=\gamma^*\|\gamma^*$ and it implies that $a^*\|q^{*\ \gamma^*\|\gamma^*}\neq0$, a contradiction. Therefore it is enough to consider the case that $a^*\|q^{*\ \gamma^*\|\gamma^*}\neq0$. Indeed, $a^*\|q^{*\ \gamma^*\|\gamma^*} \neq 0$ if and only if there exists some $b_i$ such that $b_i=\alpha$ or $b_i=\beta$ for $q=b_m\dots b_1$ and $1\leq i\leq m$. Since $\alpha$ is a source arrow and $\beta$ is a sink arrow, this is equivalent to $b_1=\alpha$ or $b_m=\beta$. If $b_1=\alpha$ (resp. $b_m=\beta$), then $a\| q$ implies $a=\alpha=q$ (resp. $a=\beta=q$) since $\alpha$ (resp. $\beta$) is a source (resp. sink) arrow. Therefore, $a^*\|q^{*\ \gamma^*\|\gamma^*} \neq 0$ if and only if $a^*\|q^*=\gamma^*\| \gamma^*$, and in this case we also have $[\gamma^*\|\gamma^*,a^*\|q^*]=0$. As a consequence,  $[\gamma^*\|\gamma^*,\mathrm{Ker}(\delta_B^1)]=[\gamma^*\|\gamma^*,Z_{spp}(\alpha,\beta)] \subset Y$, whence $[Y,\mathrm{Ker}(\delta_B^1)] \subset Y$, that is, $Y$ is a Lie ideal of $\mathrm{Ker}(\delta_B^1)$.
	
If we combine Remark $\ref{dim-Image}$ with Remark $\ref{exact-seq-of-Ker}$, then we get the following exact commutative diagram:
	\begin{center}
		\begin{tikzcd}
		& 0 & 0 \\
		& {Z_{sp}(\alpha,\beta)} & {Z_{spp}(\alpha,\beta)} & 0 \\
		0 & {Y=\psi_1(\mathrm{Im}(\delta_A^0)) \oplus Z_{sp}(\alpha,\beta)} & {\mathrm{Ker}(\delta_B^1)} & {\frac{\mathrm{Ker}(\delta_B^1)}{Y}\simeq\frac{\mathrm{HH}^1(B)}{I}} & 0 \\
		0 & {\mathrm{Im}(\delta_A^0)} & {\mathrm{Ker}(\delta_A^1)} & {\mathrm{HH}^1(A)} & 0 \\
		& {\langle \alpha \| \alpha-\beta\| \beta \rangle} & {\langle \alpha \| \alpha-\beta\| \beta \rangle} & 0 \\
		& 0 & 0
		\arrow["\iota", Rightarrow, no head, from=2-2, to=2-3]
		\arrow[from=2-3, to=2-4]
		\arrow[from=2-3, to=1-3]
		\arrow[from=2-2, to=1-2]
		\arrow[from=3-1, to=3-2]
		\arrow["{\iota_B}", from=3-2, to=3-3]
		\arrow[from=3-3, to=3-4]
		\arrow[from=3-4, to=3-5]
		\arrow[from=3-2, to=2-2]
		\arrow[from=3-3, to=2-3]
		\arrow[from=3-4, to=2-4]
		\arrow[from=4-1, to=4-2]
		\arrow["{\psi_1|_{\mathrm{Im}(\delta_A^0)}}"', from=4-2, to=3-2]
		\arrow["{\iota_A}", from=4-2, to=4-3]
		\arrow["{\psi_1}"', from=4-3, to=3-3]
		\arrow[from=4-3, to=4-4]
		\arrow["\varphi"', dashed, from=4-4, to=3-4]
		\arrow[from=4-4, to=4-5]
		\arrow[Rightarrow, no head, from=5-2, to=5-3]
		\arrow[from=5-3, to=5-4]
		\arrow[from=5-2, to=4-2]
		\arrow[from=5-3, to=4-3]
		\arrow[from=5-4, to=4-4]
		\arrow[from=6-2, to=5-2]
		\arrow[from=6-3, to=5-3]
		\end{tikzcd} $(**)$
	\end{center}
This implies that $\mathrm{HH}^1(A)\simeq \mathrm{HH}^1(B)/I$ as Lie algebras, where $I=\langle \gamma^*\|\gamma^* \rangle$ is a Lie ideal of $\mathrm{HH}^1(B)$ when we glue from the same block, and $I=0$ when we glue from different blocks.
\end{proof}

The next corollary shows that if $\mathrm{char}(k)=0$, then we can describe more precisely the relationship between the Lie structures of $\HH^1(A)$ and $\HH^1(B)$ in case of gluing a source arrow and a sink arrow from the same block of $A$. A similar result in the case of gluing a source vertex and a sink vertex can be found at the end of Corollary 3.24 in \cite{LRW}.

\begin{Cor}\label{one-dim-summand-source-sink}
	Let $A=kQ_A/I_A$ be a monomial algebra and let $B=kQ_B/I_B$ be the algebra obtained from $A$ by gluing a source arrow $\alpha:e_1 \to e_2$ and a sink arrow $\beta:e_{n-1} \to e_n$ from the same block. If $\mathrm{char}(k)=0$, then the above one-dimensional Lie ideal $I$ of $\mathrm{HH}^{1}(B)$ is contained in the center of $\mathrm{HH}^{1}(B)$ and there is Lie algebra isomorphism $$\mathrm{HH}^{1}(B)\simeq \mathrm{HH}^{1}(A)\times k.$$
\end{Cor}

\begin{proof}
	Theorem \ref{hh1-gluing-source-sink-arrows} shows that it is enough to prove the last statement. We adopt the notations in \cite[Corollary 3.24]{LRW} and the proof is parallel to it. We claim that $I$ lies in $Z(L_0)$ and $I$ is a summand of $L_0$ as Lie algebra, where $L_0=(k((Q_B)_1\|(Q_B)_1)\cap \mathrm{Ker}(\delta_B^1))/\mathrm{Im}(\delta_{(B)_0}^0)$. Then the statement follows. Indeed, if $L_0=I\oplus G$ as Lie algebras, then by the graduation of $\mathrm{HH}^{1}(B)$ one deduces that $$\mathrm{HH}^{1}(B)=L_0\oplus \bigoplus_{i\geq 1}L_i=(I\oplus G)\oplus \bigoplus_{i\geq 1}L_i=I\oplus(G\oplus \bigoplus_{i\geq 1}L_i)=:I\oplus L,$$ where $L$ is a Lie ideal of $\mathrm{HH}^{1}(B)$. Indeed,  $[I,L]=[I,G]+[I,\bigoplus_{i\geq 1}L_i]=[I,\bigoplus_{i\geq 1}L_i]\subset \bigoplus_{i\geq 1}L_i$ by \cite[Remark 2.4]{LRW} and by the above claim.  In addition, 
	\begin{equation*}
	\begin{split}
	[L,L]&=[L,G]+[L,\bigoplus_{i\geq 1}L_i]\subset [G,G]+[\bigoplus_{i\geq 1}L_i,G]+\bigoplus_{i\geq 1}L_i\\
	&\subset G+\bigoplus_{i\geq 1}L_i=L, 
	\end{split}
	\end{equation*}
	hence $[\mathrm{HH}^{1}(B),L]\subset L$. Consequently, to show $\mathrm{HH}^{1}(B)=I\oplus L$ as Lie algebras is equivalent to show $0=[I,L]=[I,\bigoplus_{i\geq 1}L_i]$.  Remark 2.4 in \cite{LRW} shows that both $\bigoplus_{i\geq 1}L_i$ and $I$ are Lie ideals of $\mathrm{HH}^{1}(B)$, hence $[I,\bigoplus_{i\geq 1}L_i]=0$. Therefore, $\mathrm{HH}^{1}(B)=I\oplus L$ as Lie algebras. Since there is an isomorphism of Lie algebras $\mathrm{HH}^{1}(A)\simeq \mathrm{HH}^{1}(B)/I$, we get $L\simeq \mathrm{HH}^{1}(A)$ as Lie algebras. 
	
Now we prove the above claim. Write $L_0=X/\mathrm{Im}(\delta_{(B)_0}^0)$, where $X=k((Q_B)_1\|(Q_B)_1)\cap \mathrm{Ker}(\delta_B^1)$ and denote by $(\bar{Q}_B)_1$ the set of equivalence classes of parallel arrows in $Q_B$. Then there is a decomposition of Lie algebras $$X=\oplus_{[a^*]\in (\bar{Q}_B)_1}X^{[a^*]},$$ where $X^{[a^*]}=\{a_i^*\|a_j^*\in \mathrm{Ker}(\delta_B^1)\mid a_i^*,a_j^*\in [a^*]\}$ is a Lie subalgebra of $\langle a_i^*\|a_j^*\in(Q_B)_1\|(Q_B)_1 \mid a_i^*,a_j^*\in [a^*]\rangle\simeq \gl_{|a^*|}(k)$. Since we glue a source arrow $\alpha:e_1\to e_2$ and a sink arrow $e_{n-1}\to e_n$, we have that $\gamma^*$ is the unique arrow from $f_1$ to $f_2$, that is, $[\gamma^*]=\{\gamma^*\}$ and $X^{[\gamma^*]}=\langle \gamma^*\|\gamma^* \rangle$. Therefore $$X=\langle \gamma^*\|\gamma^*\rangle\oplus H,$$ where $H=\oplus_{[a^*]\in (\bar{Q}_B)_1\backslash\{[\gamma^*] \}}X^{[a^*]}$. Since $[\gl_{|a^*|},\gl_{|b^*|}]=0$ for $[a^*]\neq [b^*] \in (\bar{Q}_B)_1$, we have $$[H,H]=\oplus_{[a^*],[b^*]\in (\bar{Q}_B)_1\backslash\{[\gamma^*]\}}[X^{[a^*]},X^{[b^*]}]=\oplus_{[a^*]\in (\bar{Q}_B)_1\backslash\{[\gamma^*]\}}[X^{[a^*]},X^{[a^*]}]\subset H$$ and $H$ is closed under the Lie bracket. By definition for each generator $a^*\|b^*$ of $H$,  we have $a^*\neq \gamma^*$ and $b^*\neq \gamma^*$, whence $[\gamma^*\|\gamma^*,H]=0$. It follows that $X=\langle \gamma^*\|\gamma^*\rangle\oplus H$ is also a decomposition of Lie algebras. Lemma \ref{gammapi} implies that $L_0=X/\mathrm{Im}(\delta_{(B)_0}^0)=\langle \gamma^*\|\gamma^*\rangle\oplus (H/\mathrm{Im}(\delta_{(B)_0}^0))$ and therefore $I=\langle \gamma^*\|\gamma^* \rangle$ lies in the center $Z(L_0)$ of $L_0$ and $L_0=I\oplus G$ as Lie algebras, where $G$ is the quotient of $H$ by $\mathrm{Im}(\delta_{(B)_0}^0)$. We are done.	
\end{proof}

When gluing  two arbitrary arrows $\alpha:e_1\to e_2$ and $\beta:e_{n-1}\to e_n$, the  structural results in Theorem \ref{hh1-gluing-source-sink-arrows} do not hold any more (cf. Remark \ref{dim-Image}), but we still have the following dimension formula for the first Hochschild cohomology.

\begin{Thm}\label{hh1-glue-two-arrows} {\rm(Compare to \cite[Theorem 3.20]{LRW})}
	Let $A=kQ_A/I_A$ be a monomial algebra and let $B=kQ_B/I_B$ be the algebra obtained from $A$ by gluing two arrows $\alpha:e_1 \to e_2$ and $\beta:e_{n-1} \to e_n$. Then under Assumption \ref{assum}, $$\mathrm{dim\,HH}^{1}(A)= \mathrm{dim\,HH}^{1}(B)-1-\mathrm{kspp(\alpha,\beta)}+\mathrm{sp(\alpha,\beta)}+c_A-c_B.$$ In particular, if we glue from the same block, then $\mathrm{dim\,HH}^{1}(A)= \mathrm{dim\,HH}^{1}(B)-1-\mathrm{kspp(\alpha,\beta)}+\mathrm{sp(\alpha,\beta)}$; if we glue from different blocks, then $\mathrm{dim\,HH}^{1}(A)= \mathrm{dim\,HH}^{1}(B)-\mathrm{kspp(\alpha,\beta)}$.
\end{Thm}

\begin{proof}
	This follows from Proposition \ref{Image-delta-zero-of-gluing-arrows} and Proposition \ref{Kernel-delta-one-structure}.
\end{proof}

Note that when gluing two arbitrary arrows we still have the following commutative diagram similar as in Theorem \ref{hh1-gluing-source-sink-arrows}:

\begin{center}
	\begin{tikzcd}
	& 0 & 0 & 0 \\
	0 & {Z_{sp}(\alpha,\beta)} & {Z_{spp}(\alpha,\beta)} & {\mathrm{Coker}(\varphi)} & 0 \\
	0 & {Y':=\langle\psi_1(\mathrm{Im}(\delta_A^0)),\gamma^*\|\gamma^* \rangle\oplus Z_{sp}(\alpha,\beta)} & {\mathrm{Ker}(\delta_B^1)} & {\frac{\mathrm{Ker}(\delta_B^1)}{Y'}\simeq \frac{\mathrm{HH}^1(B)}{I'}} & 0 \\
	0 & {X:=\langle \mathrm{Im}(\delta_A^0),\alpha\|\alpha \rangle} & {\mathrm{Ker}(\delta_A^1)} & {\frac{\mathrm{HH}^1(A)}{J}} & 0 \\
	0 & {\mathrm{Ker}(\psi_1|_X)} & {\langle\alpha\|\alpha-\beta\|\beta  \rangle} & {\mathrm{Ker}(\varphi)} & 0 \\
	& 0 & 0 & 0
	\arrow[from=3-1, to=3-2]
	\arrow["{\iota_B}", from=3-2, to=3-3]
	\arrow[from=3-3, to=3-4]
	\arrow[from=3-4, to=3-5]
	\arrow[from=2-1, to=2-2]
	\arrow["\iota", from=2-2, to=2-3]
	\arrow[from=2-3, to=2-4]
	\arrow[from=2-4, to=2-5]
	\arrow[from=2-4, to=1-4]
	\arrow[from=2-3, to=1-3]
	\arrow[from=2-2, to=1-2]
	\arrow["{\pi^0}", from=3-2, to=2-2]
	\arrow["{\pi^1}", from=3-3, to=2-3]
	\arrow["\pi", from=3-4, to=2-4]
	\arrow["{\iota_A}", from=4-2, to=4-3]
	\arrow[from=4-3, to=4-4]
	\arrow[from=4-4, to=4-5]
	\arrow["{\psi_1}", from=4-3, to=3-3]
	\arrow["{\psi_1|_X}", from=4-2, to=3-2]
	\arrow[from=4-1, to=4-2]
	\arrow["\varphi", dashed, from=4-4, to=3-4]
	\arrow["i", from=5-4, to=4-4]
	\arrow["\lambda", from=5-2, to=5-3]
	\arrow[from=5-1, to=5-2]
	\arrow[from=5-3, to=5-4]
	\arrow[from=5-4, to=5-5]
	\arrow[from=6-4, to=5-4]
	\arrow[from=6-3, to=5-3]
	\arrow["{i^1}", from=5-3, to=4-3]
	\arrow[from=6-2, to=5-2]
	\arrow["{i^0}", from=5-2, to=4-2]
	\end{tikzcd}
\end{center}
where $I':=Y'/\Im(\delta_B^0)$ and $J:=X/\Im(\delta_A^0)$. However, in general $Y'$ is not a Lie ideal of $\Ker(\delta_B^1)$, hence $I'$ is not a Lie ideal of $\HH^1(B)$. Also $J$ is often not a Lie ideal of $\HH^1(A)$. Therefore we can not directly compare the Lie structures of $\HH^1(A)$ and $\HH^1(B)$ using this diagram. It is also worthwhile to mention that this diagram refines to the diagram $(**)$ in Theorem \ref{hh1-gluing-source-sink-arrows}. Indeed, if $\alpha$ is a source arrow and $\beta$ is a sink arrow, then $\alpha\|\alpha \in \Im(\delta_A^0)$, which implies $X=\Im(\delta_A^0)$, so $J=0$, and  $\gamma^*\|\gamma^*\in \psi_1(\Im(\delta_A^0))$, so $Y'=Y$. Consequently, $I'=I$ and $\Ker(\psi_1|_X)=\langle \alpha\|\alpha-\beta\|\beta \rangle$, hence $\varphi$ is injective. Moreover, $\varphi$ is an isomorphism by Remark \ref{special-pair-gluing-source-sink}.

Note also that we can give another dimension formula for $\HH^1$ using the above commutative diagram. In fact,  $\dim\,\mathrm{Coker}(\varphi)$ and $\dim\,\mathrm{Ker}(\varphi)$ are equal to $\mathrm{kspp}(\alpha,\beta)-\mathrm{sp}(\alpha,\beta)$ and $1-\dim\,\mathrm{Ker}(\psi_1|_X)$ respectively, the last column of this commutative diagram gives rise to an equation as follows:
\begin{equation*}
	\begin{split}
	\dim\HH^1(A) &= \dim\HH^1(B)+\dim\, J-\dim\, I'+\dim\,\mathrm{Ker}(\varphi)-\dim\,\mathrm{Coker}(\varphi)\\ &=\dim\HH^1(B)+\dim\, J-\dim\, I'+1-\dim\,\mathrm{Ker}(\psi_1|_X)- \mathrm{kspp}(\alpha,\beta)+\mathrm{sp}(\alpha,\beta).
	\end{split}
\end{equation*}	

Although $\dim\,J$ can be described precisely, that is, $\dim\, J=0$ if $\alpha\|\alpha\in \Im(\delta_A^0)$, otherwise it is zero, we cannot give a specific description for $\dim\, I'$ or $\dim\,\mathrm{Ker}(\psi_1|_X)$, because when we glue two arbitrary arrows, $\alpha\|\alpha-\beta\|\beta$ may be not in $\Im(\delta_A^0)$ as we mentioned in Remark \ref{dim-Image} and $\dim\,\mathrm{Ker}(\psi_1|_X)$ is independent of whether $\alpha\|\alpha$ belongs to $\Im(\delta_A^0)$ or not.  In general, we have   that $\dim\, I'\in \{0,1,2,3\}$ and $\dim\,\mathrm{Ker}(\psi_1|_X)\in \{0,1\}$. However, it is interesting that $\dim\, J-\dim\, I'+1-\dim\,\mathrm{Ker}(\psi_1|_X)=c_A-c_B-1$ by comparing with Theorem \ref{hh1-glue-two-arrows}, and it is $-1$ if we glue $\alpha$ and $\beta$ from the same block, otherwise it is $0$.

\begin{Cor}\label{one-dim-summand-rad-aquare-zero} {\rm(cf. \cite[Corollary 3.28]{LRW})} 
	Let $A=kQ_A/I_A$ be a radical square zero algebra and let $B=kQ_B/I_B$ be the algebra obtained from $A$ by gluing two arrows $\alpha:e_1 \to e_2$ and $\beta:e_{n-1} \to e_n$ from the same block. If $Z_{spp}(\alpha,\beta)=0$ and if $\mathrm{char}(k)=0$, then we have Lie algebra isomorphism $$\mathrm{HH}^{1}(B)\simeq \mathrm{HH}^{1}(A)\times k.$$
\end{Cor}

\begin{proof}
	Since $Z_{spp}(\alpha,\beta)=0$, Theorem \ref{hh1-glue-two-arrows}, together with the fact that $\mathrm{sp(\alpha,\beta)}=0$ in the radical square zero case, imply that $\dim\HH^1(A)= \dim\HH^1(B)-1$. Moreover, if $\mathrm{char}(k)=0$, then Theorem 2.9 in \cite{S1} shows that $\mathrm{HH}^{1}(A) \simeq \prod_{a\in S_A} \sl_{|a|}(k)\times k^{\chi(\bar{Q}_A)}$ and $\mathrm{HH}^{1}(B) \simeq \prod_{a^*\in S_B} \sl_{|a^*|}(k)\times k^{\chi(\bar{Q}_B)}$, where $S_A$ denotes a complete set of representatives of the non-trivial classes in $(Q_A)_1$, that is, equivalence classes having at least two arrows in $(Q_A)_1$, and $\chi(\bar{Q}_A)$ denotes the first Betti number of $\bar{Q}_A$.  Note that by gluing $\alpha:e_1\to e_2$ and $\beta:e_{n-1} \to e_n$ from the same block we have $\chi(\bar{Q}_B)=\chi(\bar{Q}_A)+1$. By combining this with the condition $Z_{spp}(\alpha,\beta)=0$, we can identify the index set $S_A$ with $S_B$. Therefore we can deduce that $\mathrm{HH}^{1}(B)\simeq \mathrm{HH}^{1}(A)\times k.$
\end{proof}

\section{Center}\label{centergluingarrows}
In this section, we study the behaviour of the centers of finite dimensional monomial algebras under gluing arrows. Throughout this paper, we will denote by $Z(A)$ the center of an algebra $A$. It is well known that the 0-th Hochschild cohomology of an algebra is exactly its center, that is, there is a canonical isomorphism from $Z(A)$ to $\HH^0(A)$. Then we can identify the center $Z(A)$ with $\mathrm{Ker}(\delta^0_A)$ via the isomorphism $\HH^0(A)\simeq \Ker(\delta_A^0)$.

\begin{Def}\label{non-special-path}\rm(cf. \cite[Definition 6.1]{LRW})
	Let $A=kQ_A/I_A$ be a monomial algebra and let $B=kQ_B/I_B$ be obtained from $A$ by gluing two arrows $\alpha:e_1\to e_2$ and $\beta:e_{n-1}\to e_n$ of $A$. Let $A_{(i,j)}$ be the set of the paths from $e_j$ to $e_i$ for $1\leq i,j\leq n$ in $\mathcal{B}_A$. 
\end{Def}

\begin{Notation} We denote by:
	\begin{itemize}
		\item $\mathrm{NSp}(\alpha,\beta)$ the union of $A_{(1,n-1)},A_{(n-1,1)},A_{(2,n)}$ and $A_{(n,2)}$, 
		\item $\langle  \mathrm{NSp}(\alpha,\beta) \rangle$ the $k$-subspace of $k((Q_B)_0\| \mathcal{B}_B)$ generated by the elements $f_1 \| p^*$ and $f_2\|q^*$, where $p \in A_{(1,n-1)}\cup A_{(n-1,1)}$ and $q\in A_{(2,n)} \cup A_{(n,2)}$,
		\item $Z_{nsp}(\alpha,\beta)$ the intersection of $\langle  \mathrm{NSp}(\alpha,\beta) \rangle$ and $\mathrm{Ker}(\delta_B^0)$,
		\item $\mathrm{nsp(\alpha,\beta)}$ the dimension of $Z_{nsp}(\alpha,\beta)$.
    \end{itemize}
\end{Notation}

In general, the generators of $Z_{nsp}(\alpha,\beta)$ are $k$-linear combinations of the form $f_1 \| p^*+f_2\|q^*$ (cf. Example \ref{Eg1}, \ref{Eg2}), where $p$ is a path between $e_1$ and $e_{n-1}$ and $q$ is a path between $e_2$ and $e_n$ in $\mathcal{B}_A$, which is different from the gluing idempotents case (cf. \cite[Definition 6.1]{LRW}).

The following is a parallel version of Lemma 6.2 in \cite{LRW} in case of gluing two arrows. For some notations, we refer to Definition \ref{the-decomposition-of-delta0}.

\begin{Lem}\label{Kernel-delta-1}
	Let $A=kQ_A/I_A$ be a monomial algebra and let $B=kQ_B/I_B$ be obtained from $A$ by gluing two arrows $\alpha:e_1\to e_2$ and $\beta:e_{n-1}\to e_n$ of $A$. Then there is a decomposition as $k$-vector spaces $$\mathrm{Ker}(\delta^0_{(B)_{\geq1}})=\psi_0(\mathrm{Ker}(\delta^0_{(A)_{\geq1}})) \oplus Z_{nsp}(\alpha,\beta).$$ In particular, if we glue $\alpha$ and $\beta$ from the same block, then $\dim\Ker(\delta^0_{(B)_{\geq1}})=\dim\Ker(\delta^0_{(A)_{\geq1}})+\mathrm{nsp(\alpha,\beta)}$; if we glue $\alpha$ and $\beta$ from different blocks, then $\dim\Ker(\delta^0_{(B)_{\geq1}})=\dim\Ker(\delta^0_{(A)_{\geq1}}).$
\end{Lem}

\begin{proof}
	A direct computation shows that $\delta_B^0(\psi_0(e_i \| p))=\psi_1(\delta_A^0(e_i \| p))$ for $1\leq i\leq n$ and $$p\in (\mathcal{B}_A\backslash \{e_1,e_2,e_{n-1},e_n\}),$$ it gives rise to an injective $k$-linear map $\psi_0:\mathrm{Ker}(\delta_{(A)_{\geq 1}}^{0}) \hookrightarrow \mathrm{Ker}(\delta_{(B)_{\geq 1}}^{0})$ induced from $\psi_0: k((Q_A)_0\| \mathcal{B}_A) \to k((Q_B)_0\| \mathcal{B}_B)$. For each $\theta\in\mathrm{Ker}(\delta_{(B)_{\geq 1}}^{0})$ which lies in the complement of the subspace $\psi_0(\mathrm{Ker}(\delta_{(A)_{\geq 1}}^{0}))$, we can assume that $\theta$ is a linear combination of the elements of the form $f_1 \| p^*$ and $f_2\|q^*$, where $p$ is a path between $e_1$ and $e_{n-1}$ and $q$ is a path between $e_2$ and $e_n$ in $\mathcal{B}_A$. That is, $\theta$ is an element belongs to $\langle  \mathrm{NSp}(\alpha,\beta) \rangle$. Moreover, $\theta \in \mathrm{Ker}(\delta_B^0)$ implies that $\theta$ is an element in $ Z_{nsp}(\alpha,\beta)$. Therefore $\mathrm{Ker}(\delta^0_{(B)_{\geq1}})=\psi_0(\mathrm{Ker}(\delta^0_{(A)_{\geq1}})) \oplus Z_{nsp}(\alpha,\beta).$
\end{proof}

When we consider the relations between centers of algebras under gluing arrows, we distinguish two cases as follows: \begin{itemize}
	\item [$(i)$] when we glue from the same block, we assume that the algebra $A$ is indecomposable;
	\item [$(ii)$] when we glue from different blocks, we assume that $A$ has only two blocks, say $A_1$ and $A_2$, then we glue $\alpha:e_1\to e_2\in A_1$ and $\beta:e_{n-1}\to e_n \in A_2$.
\end{itemize}

\begin{Prop}\label{center-indec} {\rm(cf. \cite[Proposition 6.3]{LRW})}
	Let $A=kQ_A/I_A$ be an indecomposable finite dimensional monomial $k$-algebra and let $B=kQ_B/I_B$ be obtained from $A$ by gluing two arrows $\alpha:e_1\to e_2$ and $\beta:e_{n-1}\to e_n$ of $A$. Then there is an algebra monomorphism $Z(A)\hookrightarrow Z(B)$. Moreover, we have  $$\mathrm{dim}\,Z(B)=\mathrm{dim}\,Z(A)+\mathrm{nsp(\alpha,\beta)}.$$
\end{Prop}

\begin{proof}
By \cite[Remark 2.3]{LRW}, we can identify $\mathrm{Ker}(\delta_{A}^{0})$ with $Z(A)$ by $\sum e_i \| p \mapsto \sum p$ and $\sum_{i=1}^{n}e_i \| e_i \mapsto 1_A$, and similarly for $\mathrm{Ker}(\delta_{B}^{0})$ and $Z(B)$. Also notice that $\mathrm{Ker}(\delta_{A}^{0})=\mathrm{Ker}(\delta^0_{(A)_0})\oplus \mathrm{Ker}(\delta_{(A)_{\geq 1}}^{0})$ as $k$-vector spaces and the similar decomposition applies for $\mathrm{Ker}(\delta_{B}^{0})$.
	
	By Lemma \ref{Kernel-delta-1}, we know that $\psi_0$ induces an injective $k$-linear map from $\mathrm{Ker}(\delta_{(A)_{\geq 1}}^{0})$ to $\mathrm{Ker}(\delta_{(B)_{\geq 1}}^{0})$, and $\dim\Ker(\delta^0_{(B)_{\geq1}})=\dim  \Ker(\delta^0_{(A)_{\geq1}})+\mathrm{nsp(\alpha,\beta)}$. Also note that $\mathrm{Ker}(\delta^0_{(A)_0})=\langle \sum_{1\leq i\leq n}e_i \| e_i \rangle$ and $\mathrm{Ker}(\delta^0_{(B)_0})=\langle \sum_{1\leq i\leq n-2}f_i \| f_i \rangle$. We deduce that $\dim\Ker(\delta^0_{(B)_0})=\dim\Ker(\delta^0_{(A)_0})$, hence the second statement follows. Moreover, there is an injective $k$-linear map $\psi_0:\mathrm{Ker}(\delta_{A}^{0})\to \mathrm{Ker}(\delta_{B}^{0})$. Then, the fact that $p^*q^*=(pq)^*$ for $p,q \in (\mathcal{B}_A \backslash \{e_1,\cdots,e_n\})$ shows that $\psi_{0}$ gives an algebra monomorphism $Z(A)\hookrightarrow Z(B)$, and the first statement follows.
\end{proof}

 Let  $A$ be an indecomposable monomial algebra. When gluing a source vertex and a sink vertex of $Q_A$, Corollary 6.4 in \cite{LRW} gives a necessary and sufficient condition for the monomorphism between centers to be an isomorphism.
When gluing a source arrow and a sink arrow, we have the following sufficient condition for the monomorphism between centers to be an isomorphism.

\begin{Cor}\label{center-source-sink-indec}
	Let $A=kQ_A/I_A$ be an indecomposable finite dimensional monomial $k$-algebra and let $B=kQ_B/I_B$ be obtained from $A$ by gluing a source arrow $\alpha:e_1\to e_2$ and a sink arrow $\beta:e_{n-1}\to e_n$ of $A$. Then the algebra monomorphism $Z(A)\hookrightarrow Z(B)$ is an isomorphism when $A_{(n-1,2)}=\emptyset$.
\end{Cor}

\begin{proof}
	If $\alpha:e_1\to e_2$ is a source arrow and $\beta:e_{n-1}\to e_n$ is a sink arrow, then $\langle \mathrm{NSp}(\alpha,\beta)\rangle$ is generated by $f_1\|(p\alpha)^*$ and $f_2\|(\beta q)^*$, where $p,q$ are paths from $e_2$ to $e_{n-1}$ in $\mathcal{B}_A$ such that $p\alpha\neq 0$ and $\beta q\neq0$. Therefore, if there is no path from $e_2$ to $e_{n-1}$ in $\mathcal{B}_A$, that is $A_{(n-1,2)}=\emptyset$, then $Z_{nsp}(\alpha,\beta)=0$. Hence $\mathrm{Ker}(\delta^0_A) \simeq  \mathrm{Ker}(\delta^0_B)$ by Proposition \ref{center-indec}. 
\end{proof}

Note that the converse of Corollary \ref{center-source-sink-indec} is not true in general (cf. Example \ref{eg1}). If $A$ is also a radical square zero algebra, then we have the following result, which is a similar version of Corollary 6.5 in \cite{LRW}.

\begin{Cor}\label{center-rad-square-zero-indec}
	Let $A=kQ_A/I_A$ be an indecomposable radical square zero monomial $k$-algebra and let $B=kQ_B/I_B$ be obtained from $A$ by gluing two arrows $\alpha:e_1\to e_2$ and $\beta:e_{n-1}\to e_n$ of $A$. Then $Z(A)\simeq Z(B)$ if and only if there are no arrows between $e_1$ (resp. $e_2$) and $e_{n-1}$ (resp. $e_n$).
\end{Cor}

\begin{proof}
	Since $e_1,e_2,e_{n-1}$ and $e_n$ are pairwise different, the radical square zero condition yields that $Z_{nsp}(\alpha,\beta)=\langle  \mathrm{NSp}(\alpha,\beta) \rangle \cap \mathrm{Ker}(\delta_B^0)$ are generated by $f_1\|a^*$ and $f_2\|b^*$, where $a$ is an arrow between $e_1$ and $e_{n-1}$ and $b$ is an arrow between $e_2$ and $e_n$. Therefore Proposition \ref{center-indec} shows that $Z(A)\simeq Z(B)$ if and only if there are no arrows between $e_1$ (resp. $e_2$) and $e_{n-1}$ (resp. $e_n$).
\end{proof}

According to \cite{C}, the radical square zero condition implies that $\mathrm{dim}\,Z(A)=|(Q_A)_1\|(Q_A)_0|+1$. It follows from the above proof that if $A,B$ satisfy the conditions in Corollary \ref{center-rad-square-zero-indec}, then the number of loops in $Q_B$ is equal to the number of loops in $Q_A$ plus the number of arrows between $e_1$ (resp. $e_2$) and $e_{n-1}$ (resp. $e_n$).

Now we compare the centers when we glue from different blocks.

\begin{Prop}\label{center-different-block} {\rm(cf. \cite[Proposition 6.7]{LRW})}
	Let $A$ be a finite dimensional monomial algebra with two blocks $A_1$ and $A_2$. Let $B$ be a radical embedding of $A$ obtained by gluing arrows $\alpha:e_1\to e_2 \in A_1$ and $\beta:e_{n-1}\to e_n \in A_2$. Then there is an injective homomorphism of algebras $Z(B)\hookrightarrow Z(A)$. Moreover, $\mathrm{dim}\, Z(A)=\mathrm{dim}\, Z(B)+1$.
\end{Prop}

\begin{proof}
	Let $\mathcal{B}_A=\{e_1,\cdots,e_n,\alpha,\beta,a,p_1,\cdots,p_u\ |\ a\in(Q_A)_1, \mbox{the length of each }p_i\mbox{ is }\geq 2 \}$ denotes the properly chosen $k$-basis of the monomial algebra $A$. Then the subalgebra $B$ of $A$ has a $k$-basis $\mathcal{B}_B=\{e_1+e_{n-1},e_2+e_n,e_3,\cdots,e_{n-2},\gamma^*,a^*,p_1^*,\cdots,p_u^*\}$. We identify the centers $Z(A), Z(B)$ as $\mathrm{Ker}(\delta^0_A), \mathrm{Ker}(\delta^0_B)$ respectively. Let $Z(A)=Z(A)_0\oplus Z(A)_{\geq 1}$ be the decomposition corresponding to $\mathrm{Ker}(\delta_{A}^{0})=\mathrm{Ker}(\delta^0_{(A)_0})\oplus \mathrm{Ker}(\delta_{(A)_{\geq 1}}^{0})$ as $k$-vector spaces, so does for $Z(B)$.
	
	By Lemma \ref{Kernel-delta-1}, we obtain that $\mathrm{Ker}(\delta^0_{(B)_{\geq1}})\simeq \mathrm{Ker}(\delta^0_{(A)_{\geq1}})$, hence $Z(A)_{\geq 1}=\langle \sum p\ |\ \text{$p$ is a cycle in $\mathcal{B}_A$} \rangle = Z(B)_{\geq 1}$. Note that $Z(A)_0=\langle 1_{A_1},1_{A_2} \rangle$, where $1_{A_j}$ denotes the unit element in $A_j$ for $j=1,2$, and $Z(B)_0=\langle 1_B=1_{A_1}+1_{A_2} \rangle$. Therefore the canonical embedding from $Z(B)_0$ to $Z(A)_0$ induces an injective homomorphism of algebras $Z(B)\hookrightarrow Z(A)$. And the difference between dimensions of $Z(A)$ and $Z(B)$ is given by the difference between dimensions of $Z(A)_0$ and $Z(B)_0$, which is exactly one.
\end{proof}

\section{Fundamental group}
\label{fundsec}
	Let $\pi_1(Q, I)$ be a fundamental group of a bound quiver $(Q,I)$.
	Suppose that a quiver $Q$ has $n$ vertices and $m$ edges and $c$ connected components. We adopt the notation that the first Betti number of $Q$, denoted by $\chi(Q)$, equals $m - n + c$. Note that the first Betti number is equal to the dimension of the first cohomology group of the underlying graph of $Q$, see for example  \cite[Lemma 8.2]{S}. Intuitively we can say that $\chi(Q)$ counts the number of holes in $Q$.

	Recall from \cite[Lemma 1.7]{BR} that for a bound quiver $(Q,I)$ we have  $\dim\,\mathrm{Hom}(\pi_1(Q, I), k^{+}) \leq \chi(Q)$.
	Equality holds if $I$ is a monomial ideal, and more generally if $I$ is semimonomial  \cite[Section 1]{GAS} and in positive
	characteristic if $I$ is $p$-semimonomial; see after Remark 1.8 in \cite{BR}. Therefore, by Theorem C in \cite{BR}, we
	have that \[\pirank(A):=\max\{\dim\,\pi_1(Q,I)^\vee ~:~ A\simeq kQ/I, \textrm{ $I$ is an admissible ideal} \}\] is equal to $
	\chi(Q_A)$, where $\pi_1(Q,I)^\vee ~=\mathrm{Hom}(\pi_1(Q, I), k^{+})$.

The $\pirank(A)$ is a derived invariant and an invariant under stable equivalences of Morita type for selfinjective algebras, see \cite[Theorem B]{BR}. However,  it is not an invariant under stable equivalences induced by gluing idempotents \cite{LRW}.

It is interesting that the formula that compares $\pirank(A)$ and $\pirank(B)$ is the same in the case of gluing arrows case as in the case of gluing of vertices.

	\begin{Lem} \label{pi-one-rank} {\rm(cf. \cite[Lemma 5.1]{LRW})}
		Let $A=kQ_A/I_A$ be a finite dimensional monomial (or semimonomial) algebra and let $B=kQ_B/I_B$ be a finite dimensional algebra  obtained by gluing two arrows of $A$. Then $$\pirank(A)=\pirank(B)+c_A-c_B-1.$$ In particular, if we glue two arrows from different blocks, then $$\pirank(A)=\pirank(B);$$ and if we glue two arrows from the same block, then $$\pirank(A)=\pirank(B)-1.$$
	\end{Lem}

	\begin{proof} Since $A$ and $B$ are monomial algebras, we have that $\pirank(A)=\chi(Q_A)$ and $\pirank(B)=\chi(Q_B)$ by Theorem C in \cite{BR}. The statement follows from the fact the number of arrows of $Q_A$ and $Q_B$ differ by one, that is, $m_A=m_B+1$, and from the observation that $n_A=n_B+2$. Also note that the number of connected components have the relation $c_A=c_B+1$ when we glue two arrows from different blocks but $c_A=c_B$ when we glue two arrows from the same block. The same argument applies if $A$ and $B$ are semimonomial algebras.
	\end{proof}

	\begin{Rem} \rm(cf. \cite[Lemma 5.2]{LRW})
		When the characteristic of the field is positive, Lemma \ref{pi-one-rank} holds also for $p$-semimonomial algebras since the $\pirank$ coincides with the first Betti number.
	\end{Rem}

The goal of the rest of the section is to establish a deeper connection between the dual fundamental group and the first Hochschild cohomology in case of gluing a source arrow and a sink arrow. 

\begin{Lem}\label{gammapi}
Let $A=kQ_A/I_A$ be a monomial algebra and let $B=kQ_B/I_B$ be the algebra obtained from $A$ by gluing a source arrow $\alpha:e_1 \to e_2$ and a sink arrow $\beta:e_{n-1} \to e_n$ from the same block of $A$. Then  $\gamma^* \| \gamma^* \notin \mathrm{Im}(\delta_B^0)$. 
\end{Lem}

\begin{proof}
Since we are gluing from the same block of $A$, without loss of generality we can assume that $A$ is indecomposable.
    Let $B$ be the monomial algebra having monomial ideal $I_B$.
    Denote by $e_1, \dots , e_n$ the vertices of $Q_A$ and by $f_1, \dots, f_{n-2}$ the corresponding vertices of $Q_B$. Fix $f_2$ to be the base point of the fundamental group, that is, $\pi_1(Q_B, I_B)=\pi_1(Q_B, I_B, f_2)$.
    Since  $\alpha$ is a source arrow and $\beta$ is a sink arrow, there is a walk $v$ in $Q_A$ starting at $e_2$ and ending at $e_{n-1}$ that does not contain $\alpha$ or $\beta$.
    As a consequence, the corresponding walk $v^*$ in $Q_B$ will start at $f_2$ and end at $f_1$, and $v^*$ does not contain $\gamma^*$.
     Hence we can consider $\gamma^*v^*$ as an element of $\pi_1(Q_B, I_B)$. Let $g^*$ be the dual  of $\gamma^*v^*$ and let consider $g^*$ as an element of the basis of $\mathrm{Hom}(\pi_1(Q_B, I_B), k^+)$.
    
   The key idea of this proof is to show that the element $\gamma^* \| \gamma^*$ is coming from the element $g^*$ in $\mathrm{Hom}(\pi_1(Q_B, I_B), k^+)$ via the injective map defined in \cite{AdP}:
    \[
    \theta: \mathrm{Hom}(\pi_1(Q_B, I_B), k^+) \to  \mathrm{HH}^1(B).
    \]
    Recall that in order to construct $\theta$, we first need to choose for each $f_i$ a walk $w^*_{i}$ from $f_2$ to $f_i$ with $w^*_2$ being the trivial walk at $f_2$. In our case, we make the following choice: take $w^*_1= v^*$  and, for every $i\in \{3, \dots, n-2\}$, take $w^*_i$ such that they do not contain $\gamma^*$. Note that the latter choice is always possible. Indeed, if we consider a walk from $f_2$ to $f_i$ which contains $\gamma^*$, then we replace $\gamma^*$ by $(v^*)^{-1}$. Then, for each $h^*\in \mathrm{Hom}(\pi_1(Q_B, I_B), k^+)$ and for each  path  $p^*$ from $f_i$ to $f_j$, the map $\theta$ is defined as follows \cite{dPS}:
    \[
   \theta(h^*)(p^*)=h^*((w^*_j)^{-1}p^* w^*_i) p^*
    \] 
    In our case, $\theta(g^*)(\gamma^*)=g^*((w^*_2)^{-1}\gamma^* w^*_1) \gamma^*= g^*(\gamma^* v^*) \gamma^*=\gamma^*$. In addition, since each $w^*_i$ does not contain the arrow $\gamma^*$, we have that $\theta(g^*)(\delta^*)=0$ for every  arrow $\delta^*\neq \gamma^*$. Hence  $\theta(g^*)=\gamma^*||\gamma^*$.
\end{proof}

We can provide an interpretation of Theorem \ref{hh1-gluing-source-sink-arrows}, in the case of gluing a source arrow and a sink arrow, in terms of the dual fundamental group:

\begin{Prop}
    Let $A=kQ_A/I_A$ be a monomial algebra and let $B=kQ_B/I_B$ be the algebra obtained from $A$ by gluing a source arrow $\alpha:e_1 \to e_2$ and a sink arrow $\beta:e_{n-1} \to e_n$ from the same block.  Let $v$ be a walk in $Q_A$ starting at $e_2$ and ending at $e_{n-1}$ which  does not contain $\alpha$ or $\beta$. Let $v^*$ be the corresponding walk  in $Q_B$ and let $g^*$ be the dual  of $\gamma^*v^*$. Then  the following diagram is commutative:
 
\[\begin{tikzcd}
{\mathrm{Hom}(\pi_1(Q_B, I_B, f_2), k^+) }\hspace{-0.5em} && {\mathrm{Hom}(\pi_1(Q_B, I_B, f_2), k^+)/\langle g^* \rangle}\hspace{-0.5em} && {\mathrm{Hom}(\pi_1(Q_A, I_A, e_2), k^+) } \\
\\
{ \mathrm{HH}^1(B)}\hspace{-0.5em} && {\mathrm{HH}^{1}(B)/\langle \gamma^*\|\gamma^* \rangle}\hspace{-0.5em} && { \mathrm{HH}^1(A)}
\arrow["{\theta_B}", from=1-1, to=3-1]
\arrow["{\pi_2}", from=3-1, to=3-3]
\arrow["{\theta_A}"', from=1-5, to=3-5]
\arrow["\sigma"', from=1-3, to=1-5]
\arrow["\psi", from=3-3, to=3-5]
\arrow["{\pi_1}"', from=1-1, to=1-3]
\arrow["{\tilde{\theta}_B}", from=1-3, to=3-3]
\end{tikzcd}\]
where $\sigma$ and $\psi$ are isomorphisms and $\pi_1$ and $\pi_2$ are canonical projections.
\end{Prop}

\begin{proof}
In order to define $\sigma$, we first need to consider a minimal number of generators of $\pi_1(Q_A, I_A, e_2)$ and of $\pi_1(Q_B, I_B, f_2)$. 
If we consider a walk $v$ in $Q_A$ starting at $e_2$ and ending at $e_{n-1}$ which  does not contain $\alpha$ or $\beta$, then the corresponding walk $v^*$ in $Q_B$ will start at $f_2$ and end at $f_1$. Hence we can consider $\gamma^*v^*$ as an element of $\pi_1(Q_B, I_B, f_2)$. We choose the rest of the of the generators of $\pi_1(Q_B, I_B, f_2)$ as follows: let $\{v_1, \dots, v_t\}$ be a minimal number of generators in $\pi_1(Q_A, I_A, e_2)$. Since $\alpha$ and $\beta$ are source and sink arrows, respectively, such elements of the generating set will not contain $\alpha$, $\beta$ (and their formal inverses). We set $\{ v^*_1, \dots, v^*_t, \gamma^*v^*\}$ to be a minimal number of generators of $\pi_1(Q_B, I_B, f_2)$. Then we take the corresponding dual basis in $\mathrm{Hom}(\pi_1(Q_A, I_A, e_2), k^+)$ and in $\mathrm{Hom}(\pi_1(Q_B, I_B, f_2),k^+)$, respectively. More precisely, a basis of  $\mathrm{Hom}(\pi_1(Q_A, I_A, e_2), k^+)$ is given by $\{g_1,\dots,g_t\}$ where $g_i(v_i)=1$ and $g_i$ evaluated at any other element of the basis is zero. Similarly, a basis of $\mathrm{Hom}(\pi_1(Q_B, I_B, f_2), k^+)$ is given by $\{g^*_1,\dots,g^*_t, g^*\}$, where  $g^*$ is the dual  of $\gamma^*v^*$. 
The map
$\sigma$  sends $g_i^*$ to  $g_i$. 

To define $\theta_B$ we consider a very similar  parade 
data as in Lemma \ref{gammapi}, that is, take $w^*_1= v^*$, $w_2^*=f_2$,  and, 
for every $i\in \{3, \dots, n-2\}$, take $w^*_i$ such that they do 
not contain $\gamma^*$ or $(\gamma^*)^{-1}$. We define the parade 
data in $Q_A$ in terms of the parade data in $Q_B$. This choice is 
important in order to make the diagram commutative. More precisely, 
for $3\leq i\leq n-2$,  if $w^*_i$ is given by the concatenation 
of arrows and formal inverses $a^*_s\cdots a^*_1$, then  we set 
$w_i=a_s\cdots a_1$. We set $w_1=\alpha^{-1}$, $w_2=e_2$,  $w_{n-
1}=v=b_m\cdots b_1$ and  $w_n=\beta b_m\cdots b_1$, where 
$v^*=b^*_m\cdots b^*_1$. The map $\tilde{\theta}_B$ is defined by 
Lemma \ref{gammapi}. The map $\psi$ is defined in Theorem 
\ref{hh1-gluing-source-sink-arrows}. 

The left square is commutative by Lemma \ref{gammapi}.
In order to verify that the right square is commutative, take 
 $h^*$ in the basis of $ \mathrm{Hom}(\pi_1(Q_B, I_B, f_2), k^+)/\langle g^*\rangle$. For every arrow $\delta:e_i\to e_j (1\leq i,j\leq n)$ in $Q_A$, on the one hand,  
\[
{\theta}_A(\sigma (h^*))(\delta)=h(w_j^{-1} \delta w_i)\delta.
\] 
On the other hand, 
\[
\tilde{\theta}_B(h^*)(\delta^*)= h^*((w^*_j)^{-1}\delta^*w^*_i)\delta^*.
\] 
Therefore, $\psi(\tilde{\theta}_B(h^*))(\delta)=h^*((w^*_j)^{-1} \delta^* w^*_i)\delta$. Hence we should check that $h(w_j^{-1} \delta w_i)=h^*((w^*_j)^{-1}\delta^*w^*_i)$ for every arrow $\delta$. The equality holds for $\alpha$ and $\beta$ since by construction $h$ vanishes on walks that contain $\alpha$ or $\beta$. Similarly,  $h^*$ vanishes for walks that contain $\gamma^*$. Note that these will be only cases in which $w_1$, $w_n$ and their formal inverses will appear. For any other arrow $\delta:e_i\to e_j$ with $2\leq i\leq n-2$ and  $3\leq j\leq n-1$ we will have that $w_i$ corresponds to $w^*_i$ for $2\leq i\leq n-2$ and $w_{n-1}(=v)$ corresponds to $w_1^*(=v^*)$ via the map induced by $\varphi: Q_A\to Q_B$ on the set of walks. The statement follows.
\end{proof}

\section{Higher degrees}
\label{higherhoch}
In this subsection, we assume that all algebras considered are indecomposable and radical square zero.
\begin{Def}(\cite{C})
A $n$-crown is a quiver with $n$ vertices
cyclically labeled by the cyclic group of
order $n$, and $n$ arrows $a_0, \dots ,
a_{n-1}$ such that $s(a_i)=i$ and
$t(a_i)=i+1$. A $1$-crown is a loop, and a
$2$-crown is an oriented 2-cycle.
\end{Def}

Theorem 2.1 in \cite{C} provides the dimension of Hochschild cohomology: Let $Q$ be a connected quiver which is not a crown.
The dimension of the $n$-th  Hochschild cohomology group is:
$$\mathrm{dim} \HH^n(A)=|Q_n\| Q_1|-|Q_{n-1}\| Q_0|.$$

For the $n$-crown case, see Proposition 2.3 in \cite{C}. Note that there is a typo in \cite{C}  since the formula above holds for $n>1$ and not for $n>0$. On page 24 of S\'anchez Flores' PhD thesis \cite{S1}  this is corrected.
%\begin{Rem} \cite[Remark 1.5]{C}
%If the Gabriel quiver $Q$ of a radical square zero algebra $A$ has no oriented cycles, there is a path of maximum length $m$. From the above result we infer that $\HH^n(A)$ is zero for $n>m$.
%\end{Rem}

\begin{Lem} \label{parallel-path}
Let $B$ be obtained by gluing two arrows (say $\alpha$ and $\beta$) from $A$ and let $n\geq 2$. Let $\alpha_1 \in (Q_{A})_{1}$ and $p \in (Q_{A})_{n}$. If $p \| \alpha_1$, then $p^{*} \| \alpha^{*}_1$, where the map $\varphi_n:(Q_{A})_n \to (Q_{B})_n$ sends $p$ to $p^*$ (cf. Notation in Remark \ref{quiver-morphism}). In particular, $\varphi_n$ is injective. If further $\alpha$ is a source arrow and $\beta$ is a sink arrow, then $p\|\alpha_1$ implies that $\alpha_1\neq \alpha$ and $\alpha_1\neq \beta$, and $p$ does not contain $\alpha$ or $\beta$. 
\end{Lem}

\begin{proof}
The first part uses a similar argument of Proposition \ref{parallel-paths-in-monomial-algebras} (2). Since $n\geq 2$, the map
$\varphi_n$ is injective because the map $\varphi: (Q_A)_1\setminus \{\alpha,\beta\} \rightarrow (Q_B)_1\setminus \{\gamma^*\}$ sending $\delta$ to $\delta^*$ is injective, again by Proposition \ref{parallel-paths-in-monomial-algebras} (1).
We give some details for the second part of the statement. The fact that $\alpha$ is a source arrow and $\beta$ is a sink arrow implies that  $\alpha_1\neq \alpha$ and $\alpha_1\neq \beta$. In addition, $p$ does not
contain $\alpha$ or $\beta$. Let us first show that  if
$p=a_n\dots a_1$ then  $a_1\neq \alpha$. Note  that, except  $\alpha$, there is no   arrow starting at $s(\alpha)$.
 So if $a_1= \alpha$, then $\alpha_1=\alpha$ which is not possible. If $a_i=\alpha$ for $2\leq i \leq n$, then we have a
contradiction since $s(\alpha)$ is a source vertex. A
similar argument applies in the case  of the sink arrow $\beta$. 
\end{proof}

\begin{Lem}\label{parallel-paths-in-radicalsquarezero-algebras}
Let $A=kQ_A/I_A$ be an indecomposable radical square zero algebra and let $B=kQ_B/I_B$ be obtained by gluing two arrows $\alpha$ and $\beta$ of $A$. Then $\varphi_n:(Q_{A})_n \to (Q_{B})_n$ induces $k$-linear maps $\psi_{n,0}:  k( (Q_A)_n \| (Q_A)_0) \to k((Q_B)_n\|(Q_B)_0)$ and $\psi_{n,1}: k( (Q_A)_n \| (Q_A)_1) \to k( (Q_B)_n\| (Q_B)_1)$ for $n\geq 2$. In addition, $\psi_{n,1}$ is injective.
\end{Lem}

\begin{proof}
If $p\|e\in  (Q_A)_n \| (Q_A)_0$ then $s(p^*)=t(p^*)=f$, where $f$ denotes the corresponding idempotent in $B$. Hence $p^*\|f\in  (Q_B)_n \| (Q_B)_0$. By extending linearly this map we obtain the map $\psi_{n,0}$.  The construction of the map $\psi_{n,1}$ follows from Lemma \ref{parallel-path}. The injectivity of $\psi_{n,1}$ follows from the injectivity of $\varphi_n$. More precisely, if $\psi_{n,1}(p\|a)=\psi_{n,1}(q\|b)$, then $\varphi_n(p)=\varphi_n(q)$ and $a^*=b^*$. Then $a=b$ by the fact that Lemma \ref{parallel-path} gives rise to $a,b\notin \{\alpha,\beta \}$, and the injectivity of $\varphi_n$ implies that $p=q$. Therefore $p\|a=q\|b$.
\end{proof}

The injectivity of $\varphi_n$ and $\psi_{n,1}$ for $n\geq 2$ in Lemmas \ref{parallel-path} and \ref{parallel-paths-in-radicalsquarezero-algebras} holds without the assumption that  $\alpha$ is a source arrow and $\beta$ is a sink arrow. Similar results hold in \cite[Lemmas 6.9, 6.10]{LRW} in which we do not require the two idempotents to be a source and a sink.

\begin{Prop}\label{highhochgluing} {\rm(cf. \cite[Proposition 6.11]{LRW})}
Let $A$ be an indecomposable radical square zero  algebra and let $B$ be obtained from
	$A$  by gluing a source arrow  and a sink arrow
	of $A$. Then there is an injective map $\psi_n:\mathrm{Ker}(\delta^n_A) \hookrightarrow \mathrm{Ker}(\delta^n_B)$ which restricts to  $\mathrm{Im}(\delta^{n-1}_A) \hookrightarrow \mathrm{Im}(\delta^{n-1}_B)$ for $n\geq 2$ 
	(cf. Notation in  \cite[Section 2]{LRW}). In addition, $\dim\HH^n(B)-\dim\HH^n(A)\geq 0$.
\end{Prop}

\begin{proof}
Let $B$ be obtained by gluing a source arrow $\alpha$ and a sink arrow $\beta$.
Note that $\mathrm{Ker}(\delta^n_A)=k( (Q_A)_n \| (Q_A)_1) \oplus \mathrm{Ker}(D_n)$. By the proof of Theorem 2.1 in \cite{C} we know that $D_n$ is injective for $n\geq 2$ since the Gabriel quiver of $A$ is not a $n$-crown. Hence $\mathrm{Ker}(\delta^n_A)=  k( (Q_A)_n \| (Q_A)_1)$, and the map $\psi_n:\mathrm{Ker}(\delta^n_A) \rightarrow \mathrm{Ker}(\delta^n_B)$ is well defined by Lemma
\ref{parallel-path} and Lemma \ref{parallel-paths-in-radicalsquarezero-algebras}. The injectivity of $\psi_n$
follows from Lemma \ref{parallel-paths-in-radicalsquarezero-algebras}. 

By checking that $\psi_{n,1}\circ D_{n-1}=D_{n-1}\circ\psi_{n-1,0}$, we can show that $\psi_n$ restricts
to  $\mathrm{Im}(\delta^{n-1}_A) \rightarrow \mathrm{Im}(\delta^{n-1}_B)$. Let $e$ be a vertex of $Q_A$ and let $\gamma\in (Q_A)_{n-1}$ such that $\gamma$ is parallel to $e$.  On the one hand
 $$\psi_{n,1}\circ D_{n-1}(\gamma\| e)=\sum_{s(a)=e ,a\in (Q_A)_1}
 a^*\gamma^*\| a^*+(-1)^{n}\sum_{t(b)=e,b\in (Q_A)_1} \gamma^*b^*\| b^*.$$
On the other hand
$$D_{n-1}\circ\psi_{n-1,0}(\gamma\| e)= \sum_{s(a^*)=f,a^*\in (Q_B)_1}
 a^*\gamma^*\| a^*+(-1)^{n}\sum_{t(b^*)=f,b^*\in (Q_B)_1} \gamma^*b^*\| b^*.$$
Note that the vertex $e$ cannot be a source or a sink.  In  addition, we have that $e\neq t(\alpha)$ and $e\neq s(\beta)$. Indeed, except $\alpha$ there is no arrow ending at  $t(\alpha)$.  A similar argument applies to $s(\beta)$. Since for the rest of the vertices there is a bijection between the number of incoming (respectively outcoming) arrows of $Q_A$ and $Q_B$, then $\psi_{n,1}\circ D_{n-1}=D_{n-1}\circ\psi_{n-1,0}$ for every $n\geq2$.

Assume that the Gabriel quiver of $B$ is not a $n$-crown. Then by \cite[Theorem 2.1]{C} the expected inequality can be written as:
$$|(Q_B)_n\| (Q_B)_1|-|(Q_A)_n\| (Q_A)_1|\geq
|(Q_B)_{n-1}\| (Q_B)_0|
-|(Q_A)_{n-1}\| (Q_A)_0|.$$

 Let $q||f$ be an element of $k((Q_B)_{n-1}\| (Q_B)_0)$ which is not in $\mathrm{Im}(\psi_{n-1,0})$. This means that either $p=0$ or $p$ is not an oriented cycle, that is, $s(p)\neq  t(p)$ where $p^*=q\in (Q_B)_{n-1}$. Consider now $a^{*}_1q||a^{*}_1$ where $q=a^*_n\dots a^*_1$. Then $a^*_1q||a^*_1\in k((Q_B)_n\| (Q_B)_1)$ but it is not an element of  $\mathrm{Im}(\psi_{n,1})$.  In fact, if $p=0$, then $a_1p=0$. If $s(a_1)=s(p)\neq t(p)$, then $a_1p=0$. This proves the above inequality.

Assume that the Gabriel quiver of $B$ is a $n$-crown for $n\geq 1$. Then $A$ is an $A_{n+2}$-quiver.  For $A_{n+2}$, the dimensions of Hochschild cohomology groups are zero since $A_{n+2}$ is hereditary. By Proposition 2.3  and  Proposition 2.4 in \cite{C} the statement follows.
\end{proof}

Assume that we glue a source arrow and a sink arrow from the same block. The following two examples show that the difference of dimensions of higher Hochschild cohomology groups is not always one.

\begin{Ex1}
\label{exeq}
Let $A$ be radical square zero algebra with Gabriel quiver given by a zig-zag type $A_n$ quiver where  $\alpha$ is the only source arrow and $\beta$ is the only sink arrow.
% https://q.uiver.app/?q=WzAsNixbMCwwLCJlXzEiXSxbMSwwLCJlXzIiXSxbMiwwLCJlXzMiXSxbMywwLCJcXGRvdHMiXSxbNCwwLCJlX3sybi0xfSJdLFs1LDAsImVfezJufSJdLFswLDFdLFsyLDFdLFsyLDNdLFs0LDNdLFs0LDVdXQ==
\[\begin{tikzcd}
	{e_1\bullet} & {e_2\bullet} & {e_3\bullet} & \dots & {e_{2n-1}\bullet} & {e_{2n}\bullet}
	\arrow["\alpha",from=1-1, to=1-2]
	\arrow[from=1-3, to=1-2]
	\arrow[from=1-3, to=1-4]
	\arrow[from=1-5, to=1-4]
	\arrow["\beta",from=1-5, to=1-6]
\end{tikzcd}\]
Let $B$ be the radical embedding obtained by gluing the source arrow $\alpha$ and  the sink arrow $\beta$. Then $\HH^n(A)$ and  $\HH^n(B)$ are zero for $n>1$ since there are no elements of the form $|Q_n\| Q_1|$.
\end{Ex1}

\begin{Ex1}
\label{nKro}
% https://q.uiver.app/?q=WzAsNCxbMCwwLCJcXGJ1bGxldCJdLFsxLDAsIlxcYnVsbGV0Il0sWzIsMCwiXFxidWxsZXQiXSxbMywwLCJcXGJ1bGxldCJdLFswLDEsIlxcYWxwaGEiXSxbMSwyLCJcXGRlbHRhXzEiLDAseyJsYWJlbF9wb3NpdGlvbiI6NjAsIm9mZnNldCI6LTN9XSxbMSwyLCJcXHZkb3RzIiwxLHsib2Zmc2V0IjoxLCJzdHlsZSI6eyJib2R5Ijp7Im5hbWUiOiJub25lIn0sImhlYWQiOnsibmFtZSI6Im5vbmUifX19XSxbMiwzLCJcXGJldGEiXSxbMSwyLCJcXGRlbHRhX24iLDIseyJvZmZzZXQiOjV9XV0=

Let $m$ be a positive integer greater than $1$. Let $A$ be radical square zero algebra having  the following Gabriel quiver: 

\[\begin{tikzcd}
	e_1\bullet & e_2\bullet & e_3\bullet & e_4\bullet
	\arrow["\alpha", from=1-1, to=1-2]
	\arrow["{\delta_1}"{pos=0.6}, shift left=3.5, from=1-2, to=1-3]
	\arrow["\vdots"{description}, shift right=1, draw=none, from=1-2, to=1-3]
	\arrow["\beta", from=1-3, to=1-4]
	\arrow["{\delta_m}"', shift right=5, from=1-2, to=1-3]
\end{tikzcd}.\]

Let $B$ be the quiver obtained by gluing $\alpha$ and $\beta$. Let $n\geq 1$, clearly $\mathrm{dim} \HH^n(A)=0$. In addition, we have that

 \[
|(Q_B)_n\| (Q_B)_1|=
\begin{cases}
		0 & \text{if $n$ even}\\
            m^{\frac{n-1}{2}}+m^{\frac{n+3}{2}} & \text{if $n$ is odd}
\end{cases}.
  \]
Indeed, there are of no paths of even length that are parallel to an arrow in $B$. In addition, if $n$ is odd, then the number of paths of  length $n$ that are parallel to an arrow in $B$ can be counted as follows: there are $m^{\frac{n-1}{2}}$ paths of length $n$ parallel to $\gamma^*$ and there are  $m^{\frac{n+1}{2}}$ paths of length $n$ parallel to $\delta_i^*$ for $1\leq i\leq m$. Indeed,  each $\delta_i^*$ can be composed only with $\gamma^*$.
In addition, 

 \[
|(Q_B)_{n-1}\| (Q_B)_0|=
\begin{cases}
		0 & \text{if $n-1$ odd}\\
            2m^{\frac{n-1}{2}}  & \text{if $n-1$ is even} 
\end{cases}.
  \]
Indeed, there are no cycles of odd length in $B$. In addition, if  $n-1\geq 2$ is even, then the number of cycles of length $n-1$ in $B$ can be counted as follows: first note from above that  there are $m^{\frac{n-3}{2}}$ paths of length $n-2$ parallel to $\gamma^*$. So there are $m^{\frac{n-1}{2}}$ cycles parallel to $f_1$ if we compose with $\delta_i^*$ from the left side of these paths. If we compose with $\delta_i^*$ from the right side of these paths, then there are $m^{\frac{n-1}{2}}$ paths parallel to $f_2$. Therefore the above computation gives rise to ($n\ge 1$) $$\dim\HH^n(B)-\dim\HH^n(A)=|(Q_B)_n\| (Q_B)_1|-|(Q_B)_{n-1}\| (Q_B)_0|=\begin{cases}
	0 & \text{if $n$ is even}\\
	m^{\frac{n+3}{2}}-m^{\frac{n-1}{2}}  & \text{if $n$ is odd} 
	\end{cases}.$$ 
\end{Ex1}

Example \ref{nKro} shows that  although we glue a source arrow and a sink arrow from the same block, the difference between the
dimensions of $\HH^n(B)$ and $\HH^n(A)$
can be arbitrarily large for any $n>1$. In other words, for any  $n>1$ and $M>0$ we can always  find a gluing  of a source arrow  and a sink arrow from the same block such that $\dim\HH^n(B)-\dim\HH^n(A)>M$. This is very
different from the case $n = 1$, where $\mathrm{dim} \HH^1(B)-\mathrm{dim} \HH^1(A)=1$, see Theorem \ref{hh1-gluing-source-sink-arrows}.

\section{Examples}\label{Ex}

By Lemma \ref{gammapi}, if we glue a source arrow $\alpha$ and a sink arrow $\beta$ from the same block of $A$, then $\psi_1(\delta_A^0(e_1\|e_1))=\gamma^* \| \gamma^*$ is not an element in $\mathrm{Im}(\delta_B^0)$. Here is an example.

\begin{Ex1}\label{eg1}
The algebra $B$ is obtained from $A$ by gluing a source arrow $\alpha$ and a sink arrow $\beta$:
	\begin{center}
		\begin{tikzcd}
		{Q_A:} \hspace{-2em} & {e_1\bullet} & {e_2\bullet} & {e_3\bullet} & {\bullet e_4} && {Q_B:} \hspace{-2em} & {f_1\bullet} & {\bullet f_2}
		\arrow["\alpha", from=1-2, to=1-3]
		\arrow["\eta", from=1-3, to=1-4]
		\arrow["\beta", from=1-4, to=1-5]
		\arrow["{\gamma^*}", shift left=2, from=1-8, to=1-9]
		\arrow["{\eta^*}", shift left=2, from=1-9, to=1-8]
		\end{tikzcd},
	\end{center}
where $Z_A=\{\eta\alpha,\beta\eta \}$ and $Z_B=\{\eta^*\gamma^*,\gamma^*\eta^* \}$. Then, $\mathrm{Im}(\delta_B^0)=\mathrm{Im}(\delta_{(B)_0}^0)$ is generated by an element $\gamma^*\|\gamma^*-\eta^*\|\eta^*$, hence $\gamma^* \| \gamma^* \notin \mathrm{Im}(\delta_B^0)$. One can also follow the proof of Lemma \ref{gammapi}: in this case $\eta$ is the walk in $Q_A$ from $e_2$ to $e_{3}$ and $\gamma^*\eta^*$ is the walk in $Q_B$ that starts and ends at $f_2$. We denote by $g^*$ the dual of $\gamma^*\eta^*$ in $\mathrm{Hom}(\pi_1(Q_B, I_B), k^+)$. Note that $g^*$ is the only element of the basis of   $\mathrm{Hom}(\pi_1(Q_B, I_B), k^+)$. The parade data is given by $w^*_2$, which is the trivial walk at $f_2$, and $w^*_1$, which is the walk $\eta^*$. Then $\theta(g^*)(\gamma^*)=g^*(\gamma^*\eta^*) \gamma^*=\gamma^*$.  Note that  $\theta(g^*)(\eta^*)=g^*((\eta^*)^{-1}\eta^*) \eta^*$=0. Hence  $\theta(g^*)=\gamma^*||\gamma^*$.
\end{Ex1}

The following example shows that Condition (4) in Definition \ref{Special-pair} is necessary.

\begin{Ex1}\label{eg3}
The algebra $B$ is obtained from $A$ by gluing $\alpha$ and $\beta$:
	\begin{center}
		\begin{tikzcd}
		{Q_A:} \hspace{-2em} & {e_1\bullet} & {e_2\bullet} & {e_3\bullet} & {\bullet e_4} && {Q_B:} \hspace{-2em} & {f_1\bullet} & {\bullet f_2}
		\arrow["\alpha", shift left=2, from=1-2, to=1-3]
		\arrow["\eta"', from=1-4, to=1-3]
		\arrow["\beta", shift left=2, from=1-4, to=1-5]
		\arrow["{\gamma^*}", {yshift=2pt}, shift left=6, from=1-8, to=1-9]
		\arrow["{\eta^*}",{yshift=-2pt}, shift right=2, from=1-8, to=1-9]
		\arrow["a", shift left=2, from=1-3, to=1-2]
		\arrow["b"', shift right=2, from=1-4, to=1-5]
		\arrow["{a^*}",swap, {yshift=2pt},shift left=-1, from=1-9, to=1-8]
		\arrow["{b^*}", {yshift=-2pt}, shift right=6, from=1-8, to=1-9]
		\end{tikzcd},
	\end{center}
where $A$ is radical square zero and $Z_{new}=\{\eta^*a^*,b^*a^*,a^*b^* \}$. By a direct computation, we have that $$\mathrm{Ker}(\delta_A^1)=\langle \alpha\|\alpha,a\|a,\eta\|\eta,\beta\|\beta,b\|b,\beta\|b,b\|\beta \rangle,$$ $$\mathrm{Ker}(\delta_B^1)=\langle \gamma^*\|\gamma^*,a^*\|a^*,\eta^*\|\eta^*,b^*\|b^*,\gamma^*\|\eta^*,\eta^*\|\gamma^*,\gamma^*\|b^*,b^*\|\gamma^*,b^*\|\eta^*,\eta^*\|b^* \rangle$$ are 7-dimensional and 10-dimensional, respectively. Note that if we do not require the condition (4) in Definition \ref{Special-pair}, then $\mathrm{Spp}(\alpha,\beta)=\{(\alpha,\eta),(\eta,\alpha),(\beta,\eta),(\eta,\beta),(b,\eta),(\eta,b),(\alpha,b),(b,\alpha) \}$. As a result, $$\langle\mathrm{Spp}(\alpha,\beta)\rangle=\langle \gamma^*\|\eta^*,\eta^*\|\gamma^*,b^*\|\eta^*,\eta^*\|b^*,\gamma^*\|b^*,b^*\|\gamma^* \rangle,$$
\begin{equation*}
	\begin{split}
	    Z_{spp}(\alpha,\beta) & =\langle\mathrm{Spp}(\alpha,\beta)\rangle \cap \mathrm{Ker}(\delta_B^1) \\ & = \langle \gamma^*\|\eta^*,\eta^*\|\gamma^*,b^*\|\eta^*,\eta^*\|b^*,\gamma^*\|b^*,b^*\|\gamma^* \rangle.
	\end{split}
\end{equation*}
Hence $\mathrm{kspp}(\alpha,\beta)=\dim\, Z_{spp}(\alpha,\beta)=6$ yields that the dimension formula $\dim\, \mathrm{Ker}(\delta_B^1)=\dim\, \mathrm{Ker}(\delta_A^1)-1+\mathrm{kspp(\alpha,\beta)}$ in Proposition \ref{Kernel-delta-one-of-gluing-source-sink-arrows} is false. Once we include the condition (4) in Definition \ref{Special-pair}, then $\mathrm{Spp}(\alpha,\beta)=\{(\alpha,\eta),(\eta,\alpha),(\beta,\eta),(\eta,\beta),(b,\eta),(\eta,b) \}$. As a result,
$$\langle\mathrm{Spp}(\alpha,\beta)\rangle=\langle \gamma^*\|\eta^*,\eta^*\|\gamma^*,b^*\|\eta^*,\eta^*\|b^* \rangle,$$ \begin{equation*}
\begin{split}
Z_{spp}(\alpha,\beta) & =\langle\mathrm{Spp}(\alpha,\beta)\rangle \cap \mathrm{Ker}(\delta_B^1) \\ & = \langle \gamma^*\|\eta^*,\eta^*\|\gamma^*,b^*\|\eta^*,\eta^*\|b^* \rangle.
\end{split}
\end{equation*}
Therefore $\mathrm{kspp}(\alpha,\beta)=\dim\, Z_{spp}(\alpha,\beta)=4$ and the formula $\dim\Ker(\delta_B^1)=\dim\Ker(\delta_A^1)-1+\mathrm{kspp(\alpha,\beta)}$ holds.
\end{Ex1}

The next example shows that even in the case of radical square zero algebras, the number of special pairs with respect to gluing two arrows can be arbitrarily large.

\begin{Ex1}\label{eg6}
	$B$ is obtained from $A$ by gluing $\alpha$ and $\beta$:
	\begin{center}
		\begin{tikzcd}
		&            & {}                                                                                                                                    &                                                &                                                                                &             &  &      &            & {} \arrow[rd, "\ddots", phantom, shift right=0]                                                                                                                                                &    \\
		Q_A:\hspace{-2em} & e_2\bullet & \bullet e_1 \arrow[l, "\alpha"'] \arrow["a_1"', loop, distance=2em, in=125, out=55] \arrow["a_t", loop, distance=2em, in=325, out=35] & {} \arrow[lu, "\ddots", phantom, shift left=0] \hspace{-4em} & e_3\bullet \arrow["p"', loop, distance=2em, in=125, out=55] \arrow[r, "\beta"] & \bullet e_4 &  & Q_B:\hspace{-2em} & f_2\bullet & \bullet f_1 \arrow[l, "\gamma^*"'] \arrow["a_1^*"', loop, distance=2em, in=125, out=55] \arrow["a_t^*", loop, distance=2em, in=325, out=35] \arrow["p^*", loop, distance=2em, in=235, out=305] & {}
		\end{tikzcd},
	\end{center}
where $A$ is radical square zero algebra and $Z_{new}=\{p^*a_i^*,a_i^*p^* |\ 1\leq i\leq t \}$. Note that $\mathrm{Spp}(\alpha,\beta)=\{(a_i,p),(p,a_i)\ |\ 1\leq i\leq t \}$ and $\delta_B^1(a_i^*\|p^*)=\sum_{r\in Z_B}r\|r^{a_i^*\|p^*}=0$ since the length of $r^{a_i^*\|p^*}$ is equal to $2$. Similary, $\delta_B^1(p^*\|a_i^*)=0$. Consequently, $Z_{spp}(\alpha,\beta)=\langle a_i^*\|p^*, p^*\|a_i^*\ |\ 1\leq i\leq t \rangle$. Therefore $\mathrm{kspp}(\alpha,\beta)=2t$ can be arbitrarily large.
\end{Ex1}

The following example shows that in general a generator of $Z_{spp}(\alpha,\beta)$ is a $k$-linear combination.

\begin{Ex1}\label{eg2}
The algebra $B$ is obtained from $A$ by gluing $\alpha$ and $\beta$:
	\begin{center}
		\small{\begin{tikzcd}
			&                           &                                                                          &                            & \bullet e_4                                   &            &  &      &                             & \bullet f_4                                                                                                          &                                            \\
			Q_A: \hspace{-2em} & e_3\bullet \arrow[r, "b"] & \bullet e_1 \arrow[r, "\alpha"] \arrow[rr, "p", bend left, shift left=2] & \bullet e_2 \arrow[r, "c"] & \bullet e_5 \arrow[r, "\beta"] \arrow[u, "a"] & e_6\bullet &  & Q_B: \hspace{-2em} & f_3\bullet \arrow[r, "b^*"] & \bullet f_1 \arrow[r, "\gamma^*", shift left=2] \arrow[u, "a^*"] \arrow["p^*"', loop, distance=2em, in=305, out=235] & \bullet f_2 \arrow[l, "c^*", shift left=2]
			\end{tikzcd}},
	\end{center}
	where $Z_A=\emptyset$ and $Z_B=Z_{new}=\{(p^*)^2,p^*c^*,c^*\gamma^*c^*,r=a^*b^* \}$. Since both $(a,ap)$ and $(b,pb)$ belong to $\mathrm{Spp}(\alpha,\beta)$, and since $\delta_B^1(a^*\|a^*p^*)=r\|a^*p^*b^*=\delta_B^1(b^*\|p^*b^*)$, then $a^*\|a^*p^*-b^*\|p^*b^*\in Z_{spp}(\alpha,\beta)$. However, neither $a^*\|a^*p^*$ nor $b^*\|p^*b^*$ belongs to $Z_{spp}(\alpha,\beta)$.
\end{Ex1}

The next example shows that the Assumption \ref{assum} (characteristic condition) is necessary for Proposition \ref{Kernel-delta} when we glue two arbitrary arrows. It also shows that the inclusion $ Z_{sp}(\alpha,\beta)  \subseteq Z_{spp}(\alpha,\beta)$ is strict in general.
\begin{Ex1}\label{eg4}
The algebra $B$ is obtained from $A$ by gluing $\alpha$ and $\beta$:
	\begin{center}
		\begin{tikzcd}
		Q_A: \hspace{-2em} & e_1\bullet \arrow[r, "\alpha"] \arrow["\xi", loop, distance=2em, in=55, out=125] & e_2\bullet & e_3\bullet \arrow[l, "\eta"'] \arrow[r, "\beta"] & \bullet e_4 &  & Q_B: \hspace{-2em} & f_1\bullet \arrow[r, "\gamma^*", shift left=2] \arrow[r, "\eta^*"', shift right=2] \arrow["\xi^*", loop, distance=2em, in=55, out=125] & \bullet f_2
		\end{tikzcd},
	\end{center}
where $Z_A=\{\xi^2 \}$ and $Z_{new}=\{\eta^*\xi^* \}$. Since there is no special path with respect to gluing $\alpha$ and $\beta$, we have $Z_{sp}(\alpha,\beta)=0$. In this case, Assumption \ref{assum} is equivalent to the characteristic of $k$ is different from 2, and from direct computations, we get $$\mathrm{Ker}(\delta_A^1) \simeq \left\{
\begin{array}{rcl}
\langle \alpha\|\alpha,\eta\|\eta,\beta\|\beta,\xi\|\xi,\alpha\|\alpha\xi \rangle &  & {\textrm{for}\ \mathrm{char}(k)\neq 2} \\
\langle \alpha\|\alpha,\eta\|\eta,\beta\|\beta,\xi\|\xi,\alpha\|\alpha\xi,\xi\|e_1 \rangle &  & {\textrm{for}\ \mathrm{char}(k)=2}
\end{array}
\right.$$
and $$\mathrm{Ker}(\delta_B^1) \simeq \langle \gamma^*\|\gamma^*,\eta^*\|\eta^*,\xi^*\|\xi^*,\gamma^*\|\gamma^*\xi^*,\gamma^*\|\eta^*\rangle.$$ Note that when the characteristic of $k$ is equal to 2, although $\xi\|e_1\in \mathrm{Ker}(\delta_A^1)$, we have that $\psi_1(\xi\|e_1)=\xi^*\|f_1\notin
\mathrm{Ker}(\delta_B^1)$. Hence $\psi_1$ does not induce a $k$-linear map $\mathrm{Ker}(\delta_A^1)\to \mathrm{Ker}(\delta_B^1)$. Once we apply Assumption \ref{assum}, we have a canonical Lie algebra homomorphism $\mathrm{Ker}(\delta_A^1)\to \mathrm{Ker}(\delta_B^1)$ with kernel  generated by $\alpha\|\alpha-\beta\|\beta$. It is clear that $\mathrm{Spp}(\alpha,\beta)=\{(\alpha,\eta),(\eta,\alpha),(\beta,\eta),(\eta,\beta) \}$ which yields $$Z_{spp}(\alpha,\beta)=\langle \gamma^*\|\eta^* \rangle,$$ therefore $Z_{spp}(\alpha,\beta)\varsupsetneqq Z_{sp}(\alpha,\beta)$.
\end{Ex1}

Note that in case of gluing two arbitrary arrows $\alpha$ and $\beta$ under  Assumption \ref{assum}, by Proposition \ref{Kernel-delta}, there is  a Lie algebra homomorphism $\psi_1:\mathrm{Ker}(\delta^1_A) \to \mathrm{Ker}(\delta^1_B)$ with kernel generated by the element $\alpha \| \alpha - \beta \| \beta$. It is worthwhile mentioning that $\alpha \| \alpha - \beta \| \beta$ may not belong to the kernel of $\psi_1|_{\mathrm{Im}(\delta_A^0)}$. Hence we cannot deduce that $\gamma^*\|\gamma^* \in \mathrm{Im}(\delta_B^0)$, even if we glue from different blocks. The following example shows an instance of such situation:

\begin{Ex1}\label{eg5}
	The algebra $B$ is obtained from $A$ by gluing $\alpha$ and $\beta$:
	\begin{center}
		\small{\begin{tikzcd}
		&& {\bullet e_3} &&&&&&& {\bullet f_3} \\
		{Q_A:} \hspace{-2em} & {e_1\bullet} & {} & {\bullet e_2} \hspace{-3em} & {e_4\bullet} & {\bullet e_5} && {Q_B:} \hspace{-2em} & {f_1\bullet} && {\bullet f_2}
		\arrow["a", from=2-2, to=1-3]
		\arrow["\alpha", from=2-2, to=2-4]
		\arrow["b", shift left=1,from=1-3, to=2-4]
		\arrow["\beta", from=2-5, to=2-6]
		\arrow["{\gamma^*}", from=2-9, to=2-11]
		\arrow["{a^*}", from=2-9, to=1-10]
		\arrow["{b^*}", from=1-10, to=2-11]
		\end{tikzcd}},
	\end{center}
	where $A$ is radical square zero algebra, and $Z_B=\emptyset$. Consequently, we have $Z_{sp}(\alpha,\beta)=0=Z_{spp}(\alpha,\beta)$ and $$\mathrm{Im}(\delta_A^0)=\langle \alpha\|\alpha+a\|a,b\|b-a\|a,\beta\|\beta \rangle.$$ In addition, $\psi_1(\mathrm{Im}(\delta_A^0))\oplus Z_{sp}(\alpha,\beta)=\langle \mathrm{Im}(\delta_B^0),\gamma^*\|\gamma^* \rangle=\langle\gamma^*\|\gamma^*,a^*\|a^*,b^*\|b^* \rangle$ where $$\mathrm{Im}(\delta_B^0)=\langle \gamma^*\|\gamma^*+a^*\|a^*,b^*\|b^*-a^*\|a^* \rangle.$$ Therefore there is a homomorphism of Lie algebras $\psi_1|_{\mathrm{Im}(\delta_A^0)}: \mathrm{Im}(\delta_A^0) \to \psi_1(\mathrm{Im}(\delta_A^0))\oplus Z_{sp}(\alpha,\beta)$, however,  $\mathrm{Ker}(\psi_1|_{\mathrm{Im}(\delta_A^0)})=0$, where $\psi_1:\mathrm{Ker}(\delta^1_A) \to \mathrm{Ker}(\delta^1_B)$ is the Lie algebraic homomorphism having kernel generated by the element $\alpha \| \alpha - \beta \| \beta$. Also note that $\gamma^*\|\gamma^* \notin \mathrm{Im}(\delta_B^0)$ in this case which is different from the case when we glue a source arrow and a sink arrow from different blocks.
\end{Ex1}

The following two examples show that, when we glue arrows from the same block, the summands of $\delta_B^0(f_1\|p^*\gamma^*)$ and $\delta_B^0(f_2\|\gamma^*p^*)$ may cancel each other out, where $p$ belongs to $A_{(n-1,2)}$ or to $A_{(1,n)}$.

\begin{Ex1}\label{Eg1}
The algebra $B$ is obtained from $A$ by gluing a source arrow $\alpha$ and a sink arrow $\beta$:
	\begin{center}
		\begin{tikzcd}
		Q_A:\hspace{-2em} & e_1\bullet \arrow[r, "\alpha"] & \bullet e_2 \arrow[r, "\eta"] & \bullet e_3 \arrow[r, "\beta"] & \bullet e_4 & Q_B: \hspace{-2em} & f_1 \bullet \arrow[r, "\gamma^*", shift left=2] & \bullet f_2 \arrow[l, "\eta^*", shift left=2]
		\end{tikzcd},
	\end{center}
	where $Z_A=\emptyset$ and $Z_B=Z_{new}=\{\eta^*\gamma^*\eta^* \}$. Obviously, $\mathrm{NSp}(\alpha,\beta)=A_{(3,1)}\cup A_{(4,2)}=\{\eta\alpha,\beta\eta \}$ yields $\langle\mathrm{NSp}(\alpha,\beta) \rangle=\langle f_1\|\eta^*\gamma^*,f_2\|\gamma^*\eta^* \rangle$. Note that $\delta_B^0(f_1\|\eta^*\gamma^*)=\gamma^*\|\gamma^*\eta^*\gamma^*=-\delta_B^0(f_2\|\gamma^*\eta^*)$. Therefore $Z_{nsp}(\alpha,\beta)=\langle  \mathrm{NSp}(\alpha,\beta) \rangle \cap \Ker(\delta_B^0)=\langle f_1\|\eta^*\gamma^*+f_2\|\gamma^*\eta^* \rangle$.	
\end{Ex1}

\begin{Ex1}\label{Eg2}
The algebra $B$ is obtained from $A$ by gluing two arrows $\alpha$ and $\beta$:
	\begin{center}
		\begin{tikzcd}
		Q_A:\hspace{-2em} & e_1\bullet \arrow[r, "\alpha", shift left=2] & \bullet e_2 \arrow[l, "\xi", shift left=2] & \bullet e_3 \arrow[l, "a"'] \arrow[r, "\beta", shift left=2] & \bullet e_4 \arrow[l, "b", shift left=2] &  & Q_B:\hspace{-2em} & f_1\bullet \arrow[r, "\gamma^*" description] \arrow[r, "a^*" description, shift right=4] & \bullet f_2 \arrow[l, "b^*" description, shift right=8] \arrow[l, "\xi^*" description, shift right=4]
		\end{tikzcd},
	\end{center}
	where $Z_A=\{\xi\alpha\xi,b\beta b \}$ and $Z_{new}=\{a^*\xi^*,b^*a^*,\xi^*\gamma^*b^*,b^*\gamma^*\xi^* \}$. It is clear that $$A_{(1,3)}=\{\xi a,\xi ab\beta\},A_{(2,4)}=\{ab,\alpha\xi ab \}\mbox{ and }A_{(3,1)}=\emptyset=A_{(4,2)},$$ hence $\mathrm{NSp}(\alpha,\beta)=A_{(1,3)}\cup A_{(2,4)}$. Consequently, $$\langle \mathrm{NSp}(\alpha,\beta)\rangle=\langle f_1\|\xi^*a^*,f_1\|\xi^*a^*b^*\gamma^*,f_2\|a^*b^*,f_2\|\gamma^*\xi^*a^*b^* \rangle.$$ Note that the direct computations $$\delta_B^0(f_1\|\xi^*a^*b^*\gamma^*)=\gamma^*\|\gamma^*\xi^*a^*b^*\gamma^*+a^*\|a^*\xi^*a^*b^*\gamma^*-b^*\|\xi^*a^*b^*\gamma^*b^*-\xi^*\|\xi^*a^*b^*\gamma^*\xi^*=\gamma^*\|\gamma^*\xi^*a^*b^*\gamma^*,$$ $$\delta_B^0(f_2\|\gamma^*\xi^*a^*b^*)=\xi^*\|\xi^*\gamma^*\xi^*a^*b^*+b^*\|b^*\gamma^*\xi^*a^*b^*-\gamma^*\|\gamma^*\xi^*a^*b^*\gamma^*-a^*\|\gamma^*\xi^*a^*b^*a^*=-\gamma^*\|\gamma^*\xi^*a^*b^*\gamma^*$$ give rise to $Z_{nsp}(\alpha,\beta)=\langle f_1\|\xi^*a^*b^*\gamma^*+f_2\|\gamma^*\xi^*a^*b^* \rangle$.
	
\end{Ex1}

\medskip
{\bf Acknowledgements.} The first author and the third author are supported by NSFC (No. 12031014). The second author participates in the INdAM group GNSAGA.

\end{document}